\documentclass[a4paper,11pt]{amsart}
\usepackage{amsmath,amsthm,amssymb}
\usepackage{times}
\usepackage{enumerate}
\usepackage{graphicx}
\usepackage[square,sort,comma,numbers]{natbib}
\usepackage[colorlinks,linkcolor=blue,citecolor=blue,urlcolor=blue]{hyperref}
\usepackage{pstricks}
\usepackage{tikz-cd}
\usepackage{latexsym}
\usepackage{color,graphicx}
\usepackage{float}

\allowdisplaybreaks

\newtheorem{theorem}{Theorem}[section]
\newtheorem{corollary}[theorem]{Corollary}
\newtheorem{lemma}[theorem]{Lemma}
\newtheorem{proposition}[theorem]{Proposition}

\theoremstyle{definition}
\newtheorem{definition}[theorem]{Definition}
\newtheorem{remark}[theorem]{Remark}
\newtheorem{example}[theorem]{Example} 

\numberwithin{equation}{section}

\DeclareMathOperator*{\essinf}{ess\,inf}
\DeclareMathOperator{\supp}{supp}

\DeclareMathOperator{\pr}{pr}
\DeclareMathOperator{\leb}{Leb}
\DeclareMathOperator{\diag}{diag}
\DeclareMathOperator{\tr}{Tr}
\DeclareMathOperator{\Laplas}{\Delta}

\newcommand{\R}{\mathbb{R}}
\newcommand{\Rr}{\overline{\mathbb{R}}}
\newcommand{\I}{\mathbb{I}}
\newcommand{\Z}{\mathbb{Z}}
\newcommand{\N}{\mathbb{N}}
\newcommand{\F}{\mathcal{F}}
\newcommand{\p}{\mathbb{P}}
\newcommand{\B}{\mathcal{B}}

\newcommand{\Pp}{\mathcal{P}}

\newcommand{\s}{\mathcal{S}}

\newcommand{\E}{\mathbb{E}}
\newcommand{\h}{\mathcal{H}}

\newcommand{\eps}{\varepsilon}

\newcommand{\St}{\mathcal{S}^{\uparrow}}

\newcommand{\Var}{\mathrm{Var\,}}
\newcommand{\Li}{L_2^{\uparrow}}
\newcommand{\Cfb}{C^{\infty}_b}
\newcommand{\Cfo}{C^{\infty}_0}
\newcommand{\FC}{\mathcal{FC}}
\newcommand{\LL}{L_2(\Li(\xi),\Xi)}
\newcommand{\Ll}{L_2(\Xi)}
\newcommand{\D}{\mathrm{D}}

\newcommand{\e}{\mathcal{E}}
\newcommand{\Dom}{\mathbb{D}}
\newcommand{\FCC}{\mathcal{FC}_0}
\newcommand{\dm}{\mathrm{d}}
\newcommand{\cW}{\mathcal{W}}
\newcommand{\cdl}{c\`{a}dl\`{a}g }

\usepackage[colorinlistoftodos,prependcaption,color=yellow,textsize=tiny]{todonotes}

\expandafter\def\csname ver@etex.sty\endcsname{3000/12/31}

\usepackage{autonum}
\usepackage[foot]{amsaddr}
 
\usepackage{autonum}

\title[Reversible CFWD]{Reversible Coalescing-Fragmentating Wasserstein Dynamics\\ on the Real Line}

\author{Vitalii Konarovskyi$^{\dagger\ddagger\S}$} 
\author{Max-K. von Renesse$^{\ddagger}$}
\address[$\dagger$]{Fakult\"{a}t f\"{u}r Mathematik, Bielefeld Universit\"{a}t, 33615 Bielefeld, Germany}
\address[$\ddagger$]{Fakult\"{a}t f\"{u}r Mathematik und Informatik, Universit\"{a}t Leipzig, 04109 Leipzig, Germany}
\address[$\S$]{Institute of Mathematics of NAS of Ukraine, 01024 Kiev, Ukraine}
\email{vitalii.konarovskyi@math.uni-bielefeld.de}
\email{renesse@uni-leipzig.de}
\date{\today}



\begin{document}

\maketitle

\begin{abstract}
We introduce a family of reversible fragmentating-coagulating processes of 
particles of varying size-scaled diffusivity with strictly local interaction  
on the real line as mathematically rigorous description of colloidal motion of fluids. The associated measure valued process provides a weak solution to a corrected Dean-Kawasaki equation for supercooled liquids without dissipation. Our construction is based on the introduction and analysis of a fundamentally new family of equilibrium measures for the associated dynamics and their Dirichlet forms. We identify the intrinsic metric as the quadratic Wasserstein distance, which makes the process a non-trivial example of Wasserstein diffusion.
\end{abstract}

\section{Introduction and Statement of main results}
\subsection{Motivation}
This paper is a continuation in a series of studies started in 
\cite{MR2537551} when we asked for natural generalizations of Brownian motion 
of a single point to the case  of an infinite or diffuse interacting 
particle 
system with conserved total mass. As critical  consistency condition with 
respect to  the 
trivial  case of the empirical (Dirac) measure following a single Brownian 
motion  we 
put the requirement that the local fluctuations of any such 
probability measure valued diffusion $\{\mu_t\}_{t \geq 0}\in \mathcal P(\R^d)$ 
be governed by 
a Varadhan formula of the form 
\[ \p \{\mu_{t+\eps} \in A\} \sim \exp\left ({- \frac{d^2_{\mathcal W} 
(\mu_t, A )}{2\eps} 
}\right), \quad \eps \ll 1, \quad A \subset \mathcal P(\R^d),\]    
where $d_{\mathcal W }$ denotes the quadratic Wasserstein distance on 
$\mathcal 
P (\R^d)$.

 Physically, this means that the spatial 
fluctuations of 
such a measure valued process  $\mu_\cdot$  should become high at locations 
where density of 
$\mu_t$ is low and vice versa, i.e. scaling of diffusivity is inverse 
proportional to density. On the level of mathematical heuristics we can combine 
the required Wasserstein Varadhan formula with  Otto's formal infinite 
dimensional Riemannian picture of optimal transport \cite{Felix:2001} to obtain 
SPDE models of the form   \[ d\mu_t = F(\mu_t)dt + \mbox{div}( \sqrt \mu_t 
dW_t), \quad \mu_t \in 
\mathcal P (\R^d), \]
where  $dW_{\cdot}$ is a white noise vector field on $\R^d$ and $F$ is a model 
dependent drift operator. The canonical choice 
\[ F(\mu_t) = \beta \Delta \mu_t, \quad \beta \geq 0, \] 
yields  the so called \textit{Dean-Kawasaki equation} for supercooled liquids 
appearing in the physics literature \cite{Dean:1996,Kawasaki199435,MR978701,RevModPhys.87.593,1742-5468-2014-4-P04004,doi:10.1063/1.4883520,1742-5468-2016-11-113202,doi:10.1063/1.478705,Spohn:1991} (see also~\cite{Zimmer:2018,Gess:2017,Cornalba:2021,Cornalba:arxiv:2022,Gess_SMFE:2022,Fehrman:arxiv:2021,Rotskoff:arxiv:2018,Wu_Zhang:arxiv:2022} for the regularised versions of the Dean-Kawasaki equation and~\cite{Helfmann:2021,Torre:2015,Kim:2017,Cornalba:arx:2021} for the numerical investigation) but in~\cite{Konarovskyi_DK:2018,Konarovskyi_DK_smooth_drift:2018} we show that this equation is either trivial or ill posed, depending on the value of $\beta$.  
However,  as shown in \cite{MR2537551,Renesse:2010}, in $d=1$ for 
$\beta >0$, and more 
recently in \cite{Konarovskyi_LDP:2015} for $\beta =0$, the model has 
non-trivial  martingale 
solutions if one 
admits a certain additional nonlinear drift operator $ 
\Gamma_\beta(\mu_t)dt $ as 
correction. 
The 
correction is the same for all $\beta >0$ such that we 
arrive at the family of models 
\[ d\mu_t= \beta \Delta \mu_t dt + \Gamma_{i} (\mu_t) dt + \mbox{div}( \sqrt \mu_t 
dW_t),\]
where $i \in \{0,1\}$ depending whether $\beta =0 $ or $\beta >0 $. The two 
expressions for $\Gamma_0$ and $\Gamma_1$ are similar, but  
the constructions of the solutions for the two cases are 
very different. In \cite{MR2537551}
 we use abstract Dirichlet form methods, in 
\cite{Konarovskyi_LDP:2015}
 we construct an explicit system of a continuum of coalescing 
Brownian  
particles of infinitesimal initial mass which slow (i.e.\ cool) down as they 
aggregate to bigger and bigger macro-particles
before they eventually collapse to a single Brownian motion. At positive time 
the system consists of finitely many particles of different sizes almost surely, 
such that the distribution 
\[ \Gamma_0(\mu_t) = \frac{1}{2}\sum_{z \in \supp(\mu_t)} (\delta_z)''  \]
is well defined for $t>0$.

The point of departure of this work is the question whether there is a 
reversible counterpart to the coalescing particle model for the $\beta 
=0$ case. In terms of the analogy to the Arratia flow \cite{Arratia:1979} (see also \cite{MR3433579,MR3055261,MR2329772,MR2671379,Riabov:2017,MR2840266,MR2259212,MR1785393,MR2113855,MR1955353,MR2068474,Le_Jan:2004,MR2052863,MR2094432,MR2408586,MR31557821})  
this 
means that we ask 
for a Brownian Net \cite{MR2408586} type extension of the modified massive Arratia 
flow from \cite{Konarovskyi:2014:arx,Marx:2018,Konarovskyi_LDP:2015,Konarovskyi:CD:2020} which  should
then include also particle break-ups but  still  satisfies the 
characteristic 
scaling requirement regarding the diffusivity of the aggregate particles. We note that a particle model without interaction in dimension $d\geq 2$ which satisfies a similar martingale problem was considered in~\cite{Schiavo2018}.

\subsection{Heuristic Description of the Model}
 The 
main result of this work is an affirmative answer. We give it by constructing 
in rather explicit way  
a new family of measure valued processes 
on the real line which solve the 
same martingale problem for $\beta =0$ and $\Gamma_i=\Gamma_0$ as the modified massive  
Arratia flow in \cite{Konarovskyi_LDP:2015}, which satisfy the 
Wasserstein Varadhan formula and which are reversible. In this sense the new 
processes interpolate between the two 
previously known models. 

As in the case of the  modified Arratia flow, the  model describes the 
motion of an  uncountable collection of particles which are parametrized by the 
unit interval as index set and  move  on the 
real axis. It is 
assumed that the initial parametrization is monotone in particle location. The 
dynamics  
 will preserve the monotone alignment,  hence a state 
of 
the system at time $t$ is given by a monotone real function $X_t: 
(0,1)\mapsto \R$, i.e. $X_t(u)$ is the position of particle $u$ at time $t$. 
The corresponding empirical measure of the state is given by $\mu_t := 
(X_t)_{\#}(\leb)\in \mathcal P(\R)$  (image measure  of Lebesgue measure $\leb$ on $[0,1]$ under $X_t$). We  call the  atoms of $\mu_t$  empirical particles, the size 
of an atom located in $x \in \R$ at time $t$  given by $m(x,t) = \leb\{u \in 
(0,1):\ X_t(u) = x \}$. 

The basic idea for the construction of $\mu_{\cdot}$ is to 
use (sticky) reflection interaction when particles are at the same location. As 
for the 
'stickiness', particles sitting at the same location will be subject to the 
same 
random, i.e. Gaussian perturbation of their location. Since they share a 
common perturbation the net volatility of this 
perturbation is scaled in inverse proportional way by the total mass of 
particles occupying the same spot, i.e. the size of the empirical particle at 
that location.  Second, the random perturbations at 
different spots are independent. 

For the 'reflection' part of the interaction we assign once and for all 
times to each particle a certain number 
\[ [0,1] \ni u \mapsto \xi(u) \in \R,\] 
which we call its \textit{interaction 
potential}. The function $\xi$ is a free parameter of the model. 

 In addition to the random 
forcing described above, each particle will also experience a drift force given 
by the  
difference between its own interaction potential and the average interaction 
potential among all particles occupying the same location. As a consequence, if 
all occupants of a certain spot have the same interaction potential, none of 
them will feel any drift. (As they also share 
the same random forcing, in this case they will move but  stay together for all 
future 
times.) 
Conversely, big 
differences in interaction potential lead to strong drift apart among the 
particles sitting at the same location. 

The most physical  choice for $\xi$ is that of a linear function 
$\xi(u) = \lambda u$ with some $\lambda \geq 0$. In this case the break-up  
mechanism for an empirical particle depends 
only on its size. As a 
result, 
$\lambda$ controls the strength of the break-up mechanism.

Below is a  simulation of the empirical measure  process $\mu_t$, $t\geq 0$, 
for $\xi = \mbox{id}$ starting from $\mu_0=\delta_0$. Grayscale colour coding 
is for atom sizes. The red line is  the 
center of mass of the system which is always a standard Brownian motion regardless the 
choice of $\xi$.
\begin{center}
\includegraphics[height=7cm, width=11.5cm]{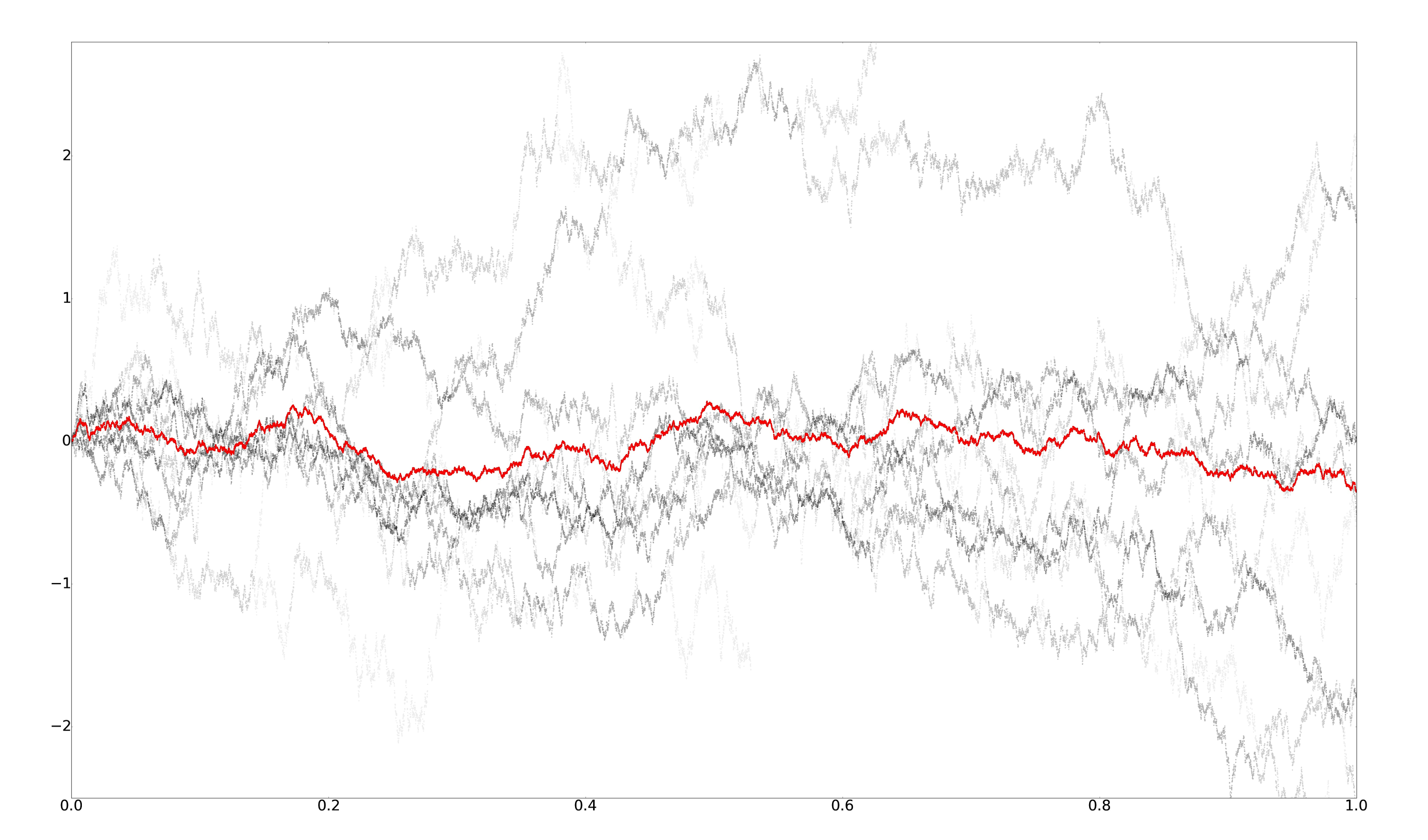}
\end{center}

\subsection{Rigorous statement of main results}

We will present now our main result in a rigorous fashion  
in terms 
of the measure valued process $\mu_{\cdot}$ assuming values in the set 
$ \Pp_2(\R)$ of  Borel probability measures on the real 
line with finite second moment and equip it with $2$-Wasserstein distance $d_{\cW}$ (see~\eqref{equ_wasserstein_distance} for the precise definition of $d_{\cW}$). 

The free parameter of the model is 
given in terms of some $\eta \in \Pp_2(\R)$, or equivalently by the choice of 
$\xi= g_\eta$, where for $\rho \in \Pp_2(\R)$ we denote by $g_\rho$ its right 
continuous quantile function, i.e 
\[ [0,1] \ni u \mapsto g_\rho (u) := \inf\{ x \in \R: \rho((-\infty, x]) > 
u\}.\] 

Given $\eta \in \Pp_2(\R)$ we introduce 
the set of all monotone transformations of $\eta$, i.e.
\[ \Pp_2^\eta(\R) := \{ \rho \in \Pp_2(\R) :\ \rho = h_{\#}(\eta) \mbox{ for 
some non decreasing } h: \R\mapsto \R\}, \]
which is a $w_2$-closed subset of $\Pp_2(\R)$. Finally, we write  
\[\Pp_2^a(\R) = \left\{  \rho=\sum_{k=1}^{n} a_k 
\delta_{z_k}:\  \sum_{k=1}^{n} a_k =1,\  a_k>0,\ z_k\in \R,\ k=1,\dots,n,\ n \in \N\right\}  \] 
for  the subset of purely countably atomic probability measures on $\R$, 
and for $\rho \in \Pp_2^a (\R)$ we set  
\[ |\rho| = \sum_{z\in\supp\rho} \delta_{z} \in \Pp_2 ( \R).\]

Below we will work with the algebra of ('smooth') functions $\FC$ on 
$\Pp_2(\R)$ which is generated by functions of the form 
\begin{align} F(\rho) &= u\left(\langle g_\rho,h_1\rangle,\ldots,\langle g_\rho,h_m\rangle\right) \cdot  \varphi (\|g_\rho\|_2^2)
\end{align}
where $u\in\Cfb(\R^m)$, $\varphi\in\Cfo(\R)$, $h_i \in 
L_2:=L_2[0,1]$, $i\in[m]$, $\langle \cdot,\cdot\rangle$ denotes the 
standard $L_2$-inner product and $\|\cdot\|_2$ is the norm on $L_2$. Writing $F(\rho) = \Phi(g_\rho)$ for $F\in \FC$,  we define the \textit{gradient} of $F \in \FC$ by  
\[ \D F(\rho):= \pr_{g_\rho} \nabla^{L_2} \Phi(g_\rho),\] 
where $\nabla^{L_2} \Phi$ denotes the standard $L_2$-gradient of $F$ which is 
defined by 
\[ \langle \nabla^{L_2}  \Phi(g) , h \rangle = \frac{\partial}{\partial \eps} 
\Phi(g+\eps h )|_{\eps=0}, \quad \forall h \in L_2, \]
and $\pr_{g}$ denotes the orthogonal projection in $L_2$ onto the subspace 
of functions which are measurable with respect to the $\sigma$-algebra 
$\sigma(g)$ on $[0,1]$ generated by the function $g$. We will 
also use  the projection $\pr^\bot$ to the complement, i.e. $\pr_g^\bot h 
= h -\pr_gh$. We will also denote the integration of a function $\psi$ with respect to a measure $\rho$ by $\langle \psi,\rho\rangle$.

With these preparations we can 
summarize the main result of this paper as follows. 
\begin{theorem} For $\eta \in \Pp_2(\R)$ with compact support there exists a measure $\Xi^\eta$ on 
$\Pp^2(\R)$  
with $\supp\Xi^\eta = \Pp_2^\eta(\R)$ such that the quadratic 
form 
\[ \e(F,F) = \int_{\Pp_2^\eta(\R)} \|\D F(\rho) \|_{2}^2 \, 
\Xi^\eta(d\rho), \quad F \in \FC, \] 
is closable on $L_2(\Pp_2^\eta, \Xi^\eta)$, its closure being a local 
quasi-regular Dirichlet form on $L_2(\Pp_2^\eta, \Xi^\eta)$. 

Let $\mu_t$, $t \in [0, \zeta)$, the properly associated 
$\Pp_2^\eta(\R)$-symmetric diffusion process with life time $\zeta>0$. Then 
\begin{enumerate}
 \item[i)] for almost all $t \in [0,\zeta)$ it holds that $\mu_t\in \Pp_2^a$ almost 
surely;
\item[ii)] for all $f \in \Cfo(\R)$ the process 
\[ M^f := \langle \mu_t , f\rangle - \frac{1}{2}\int_0^t \langle  
f'',|\mu_s|\rangle ds  \] 
is a local martingale with finite quadratic variation process 
\[ [M^f]_t = \int_0^t \langle (f')^2,\mu_s \rangle ds;\]
\item[iii)] for all $h \in L_2$ the process 
\[ \tilde M^h :=  \langle g_{\mu_t} , h\rangle  -\frac{1}{2} \int_0^t \langle 
\pr_{g_{\mu_s}}^\bot h, g_\eta \rangle ds  \] 
 is a local martingale with finite 
quadratic variation process 
\[ [\tilde M^h]_t = \int_0^t \|\pr_{g_{\mu_s}}\!\! h\|_2^2 ds;\]

\item[iv)] for all measurable $A, B \subset \Pp_2^\eta$ with $0 < \Xi^\eta(A)  
\Xi^\eta(B) < \infty$ and $A$ or $B$ open it holds that 
\[ \lim_{t \to 0}t \cdot \ln \p(\mu_0 \in A, \mu_t \in B ) = - 
\frac{d_{\mathcal W }^2 (A,B) }2,\]
where $d_{\mathcal W } (A,B) = \essinf_{(\rho, \lambda) \in A\times B} 
d_{\mathcal W}(\rho, \lambda)$. 
\end{enumerate}
\end{theorem} 

\begin{remark}
\begin{enumerate}
    \item[1)] Property ii) in the theorem  
above is equivalent to saying that $\mu_{\cdot}$ is a martingale solution to the 
SPDE 
\[d\mu_t = \Gamma_0(\mu_t) dt + {\rm  div}\, (\sqrt \mu_t dW_t)\] 
if one works with the canonical set of test functions of the type $\rho \mapsto 
\Phi(\rho) := \varphi(\langle f, \rho\rangle )$ with $\varphi, f \in 
C_0^\infty(\R)$. This collection of test functions  is commonly used in 
the theory of measure valued diffusion processes.  
Since ii) holds true regardless the choice of $\eta \in \Pp_2(\R)$, it is 
clearly not sufficient to characterize the process $\mu_{\cdot}$. This shows in 
particular that the martingale problem encoded by ii) alone 
is not well-posed. For instance,  the  
solution given by the modified  Arratia flow in \cite{Konarovskyi_LDP:2015}
 is obtained by choosing $\eta 
= \delta_z$ for some $z \in \R$, which, however, is not reversible. 

\item[2)] In fact, property ii) will be a rather straightforward  consequence of the 
stronger assertion iii),  which is 
equivalent to the statement that process $X_t := g_{\mu_t}$, $t \in [0, \zeta)$, 
is a weak solution to the SDE in infinite dimensions 
\begin{equation}\label{f_SDE_for_X_intro}
    dX_t = \frac{1}{2}\pr_{X_t}^\bot \xi \,dt + \pr_{X_t} dW_t,
\end{equation}
where $\xi = g_\eta$ and $dW$ is $L_2$-white noise. This representation is the 
justification for the heuristic description of the model in the 
previous section. As discussed in \cite{Konarovskyi_LDP:2015} the modified 
massive Arratia flow solves the same SDE with $\xi ={\rm const}$, i.e. $\eta= 
\delta_z$ for some $z \in \R$. 

\item[3)] Property iii) together with the fact that $\supp{\Xi^\eta 
}=\Pp_2^\eta$ imply in particular that the process $\mu_{\cdot}$ explores the entire 
$\Pp_2^\eta$-space. Note that $\Pp_2^\eta = \Pp_2$ iff  $\eta$ has no atoms.

\item[4)] In Section~\ref{section_quasi_regularity}, we give a first condition assuring infinite lifetime 
$\zeta = \infty$. This will be the case if e.g. $\eta ([a,b])=1$ for some 
$a\leq b$ and $\eta(\{a\}) \cdot \eta(\{b\})>0$. 
\end{enumerate}
\end{remark}

\begin{remark} Our construction given in the subsequent sections is strongly 
related to diffusion processes on domains with so called sticky-reflecting 
boundary conditions. In fact, as in \cite{MR2537551} we will 
cast the  measure valued process 
$\mu_{\cdot}$ in terms of the associated process of quantile functions $X_{\cdot}= 
g_{\mu_{\cdot}}$,  assuming values in the set $D^{\uparrow}$  of 
non decreasing functions on 
$[0,1]$. We view $D^{\uparrow}$ as a closed convex cone embedded 
in the topological space 
$L_2$. As our main and critical step  we construct the measure 
$\Xi=\Xi^\xi$ on $D^\uparrow$ which allows for an integration by 
parts formula to obtain a closable pre-Dirichlet form 
\[ \e(F,F) = \int_{D^\uparrow} \|\D F(g)\|_{L_2}^2 \Xi(dg).\]

As a subset of $L_2$ the space  $D^{\uparrow}$ has no interior since  
$\partial D^{\uparrow}$ is dense in $D^{\uparrow}$, hence we need a 
non-standard construction of a candidate measure $\Xi$. Our approach is to 
define $\Xi$ on the subset $\St$ 
of piecewise constant non decreasing functions. The set $\St= 
\bigcup_{n=0}^{\infty} \St_n$ has a natural structure as a generalized  non locally 
finite 
simplicial complex, where each $\St_n$ is the collection of all piecewise 
constant $n$-step functions. In this picture each connected component of the 
relative affine interior of 
$\St_n$ can be viewed as an $n$-dimensional face of $\St$ which is the common 
boundary of uncountably many $(n+1)$-dimensional faces that are parametrized by 
points in appropriate simplex. The measure $\Xi^\xi$ is then obtained 
by putting an $n$-dimensonal measure  $\Xi^\xi_n$ on each $\St_n$  for all $n$ 
in a way which is consistent with the hierarchical structure of $\St$. As 
a result we obtain a measure on a simplicial complex with positive mass on all 
faces of arbitrary dimension. In this picture  the gradient operator appearing 
in the Dirichlet form above is obtained as projection of the full gradient to 
the effective tangent space on the respective faces and is therefore 
geometrically natural. The outcome is a  Dirichlet form which generalizes the 
case considered e.g.  in \cite{Grothaus:2014} to the  (infinite dimensional) 
case of sticky-reflecting behaviour in piecewise smooth domains along embedded 
boundaries but now of arbitrary codimension.  
\end{remark}

The \textit{structure of this work} is as follows. After some preliminaries we 
start off in Chapter~\ref{section_finite_system} by introducing the  model in a special 
case when the system consists of a fixed finite number of atoms with prescribed 
masses. The atoms can coalesce and fragmentate, but fragmentation is allowed 
only in accordance with the initially assigned mass portions. This chapter 
exhibits the basic mechanism of the system in a finite dimensional 
situation. Section~\ref{section_measure} contains the construction of the measure $\Xi^\xi$ in the general 
case. We identify its support and show certain moment bounds which are critical 
for the quasi-regularity of the Dirichlet form which we introduce in Section~\ref{section_dirichlet_form}. The core result of Section~\ref{section_dirichlet_form} is the integration by parts formula which is needed for closability. In Section~\ref{section_quasi_regularity} we establish quasi-regularity. We also show conservativeness in a special case. Section~\ref{intrinsic metric} is devoted to the 
identification of the intrinsic metric which leads to the desired Varadhan 
formula by applying a general theorem by Ariyoshi and Hino 
\cite{Ariyoshi:2005}. In Section~\ref{section_process} we wrap up the results in terms of the induced measure valued process and the related martingale problem.

\section{Preliminaries}\label{sec_notation}

For $p\in[1,\infty]$ we denote the space of all $p$-integrable (essentially 
bounded if $p=\infty$) functions (more precisely equivalence classes) from 
$[0,1]$ to $\R$ with respect to the Lebesgue measure $\leb$ on $[0,1]$ by $L_p$ 
and $\|\cdot\|_p$ is the usual norm on $L_p$. The inner product in $L_2$ is 
denoted by $\langle\cdot,\cdot\rangle$. Let $D^{\uparrow}$ be the set of 
c\`{a}dl\`{a}g non decreasing functions from $[0,1]$ into 
$\Rr=\R\cup\{-\infty,+\infty\}$. For convenience, we assume that all functions 
from $D^{\uparrow}$ are continuous at $1$. Let $L_p^{\uparrow}$ be the subset 
of 
$L_p$ that contains functions (their equivalence classes) from $D^{\uparrow}$.

Note that $L_2^{\uparrow}$ is a closed subset of $L_2$, by~\cite[Corollary~A.2]{Konarovskyi:2017:EJP}. Consequently, $L_2^{\uparrow}$ is a 
Polish space with respect to the distance induced by $\|\cdot\|_2$.
 
If $f=g$ a.e., then we say that $f$ is a \textit{modification} or 
\textit{version} of $g$ or $g$ is a \textit{modification} or \textit{version} 
of 
$f$.

\begin{remark}
Since each function $f$ from $L_2^{\uparrow}$ has a unique modification from 
$D^{\uparrow}$ (see, e.g.,~\cite[Remark~A.6]{Konarovskyi:2017:EJP}), 
considering 
$f$ as a map from $[0,1]$ to $\Rr$, we always take its modification from 
$D^{\uparrow}$.
\end{remark}

We set for each $n\in\N$
$$
E^n=\{x=(x_1,\ldots,x_n)\in\R^n:\ x_i\leq x_{i+1},\ i\in[n-1]\}
$$
and
$$
E_0^n=\{x=(x_1,\ldots,x_n)\in\R^n:\ x_i< x_{i+1},\ i\in[n-1]\},
$$
where $[n]=\{1,\ldots,n\}$. Also let
$$
Q^n=\{q=(q_1,\ldots,q_{n-1}):\ 0<q_1<\ldots<q_{n-1}<1\}
$$
for all $n\geq 2$. Considering $q$ from $Q^n$, we will additionally take 
$q_0=0$ and $q_n=1$.

Next, for $g\in\Li$ we denote the number of distinct values of the function $g\in D^{\uparrow}$ by $\sharp g$. 
If $\sharp g<\infty$, then $g$ is called a \textit{step function}. The set of all step functions is denoted by $\St$. 
\begin{remark}\label{rem_repres_of_step_functions}
If $\sharp g=n$, then there exist unique $q\in Q^n$ and $x\in E^n_0$ such that
$$
g=\sum_{i=1}^nx_i\I_{[q_{i-1},q_i)}+x_n\I_{\{1\}},
$$
where $\I_A$ is the indicator function of a set $A$.
\end{remark}

If $E$ is a topological space, then the Borel $\sigma$-algebra on $E$ is 
denoted 
by $\B(E)$.

For any family of sets $\h$ we denote the smallest $\sigma$-algebra that 
contains $\h$ by $\sigma(\h)$. Similarly, $\sigma(f)=\sigma(\{f^{-1}(A):\ 
A\in\B(\R)\})=\{f^{-1}(A):\ A\in\B(\R)\}$
for a function $f$ taking values in $\R$. For $g\in \Li$ let 
$\sigma^{\star}(g)$ 
denote the completion of the $\sigma$-algebra $\sigma(g)$ with respect to the 
Lebesgue measure on $[0,1]$ and $\pr_g$ be the orthogonal projection operator in $L_2$ on 
the closed linear subspace 
$$
L_2(g):=\{f\in L_2:\ f\ \mbox{is}\ \sigma^{\star}(g)\mbox{-measurable}\}.
$$
By~\cite[Lemma~1.25]{Kallenberg:2002}, $\sigma^{\star}(g)$ and $L_2(g)$ are 
well-defined for each equivalence class $g$ from $\Li$. Also we set 
$\Li(g)=L_2(g)\cap\Li$.

\begin{remark}\label{rem_proj_and_ecpect}
\begin{enumerate}
 \item[(i)]  For each $h\in L_2$ the function $\pr_gh$ coincides with the 
conditional expectation $\E(h|\sigma^{\star}(g))$ on the probability space 
$([0,1],\mathcal{L}([0,1]),\leb)$, where $\mathcal{L}([0,1])$ denotes the 
$\sigma$-algebra of Lebesgue measurable subsets of $[0,1]$. 
 
 \item[(ii)] For each $h\in L_2$, $\E(h|\sigma^{\star}(g))=\E(h|\sigma(g))$ a.e.
 
 \item[(iii)] The projection $\pr_g$ maps the space $\Li$ into $\Li$, by~\cite[Lemma~A.4]{Konarovskyi:CFWDPA:2022}.
\end{enumerate}
\end{remark}

\section{Finite system of sticky reflected diffusion particles}\label{section_finite_system}

The aim of this section is to construct a finite system of diffusion particles on the real line with sticky-reflecting interaction. Also this section gives a motivation for the definition of the system in the general case. We will use a Dirichlet form approach. In particular, we use ideas from~\cite{Grothaus:2014} for the description of the sticky-reflecting mechanism. Let $n\in\N$ and $m_i\in(0,1]$, $i\in[n]$, with $m_1+\ldots+m_n=1$ be fixed. That numbers will play a role of a number of particles and particle masses, respectively.

\subsection{Some notation}

Let $\varTheta^n$ denote the set of all ordered partitions of $[n]$. We take $\theta=(\theta_1,\ldots,\theta_p)\in\varTheta^n$ and denote the number of sets in the partition $\theta$ by $|\theta|$, i.e. $|\theta|=p$.  Let
$$
E_{\theta}=\{x\in E^n:\ x_i=x_j\ \ \Leftrightarrow\ \ i,j\in\theta_k\ \ \mbox{for some}\ \ k\in[|\theta|]\}.
$$
Remark that the sets $E_{\theta}$, $E_{\theta'}$ are disjoint for $\theta\neq\theta'$ and $E^n=\bigcup_{\theta\in\varTheta^n}E_{\theta}$.

Let $R_{\theta}$ be the bijection between $E_{\theta}$ and $E^{|\theta|}$ defined as follows
$$
R_{\theta}\left(x_1,\ldots,x_n\right)=(y_1,\ldots,y_{|\theta|}),
$$
where $y_k=x_i$ for some $i\in\theta_k$ (and, consequently, for all $i\in\theta_k$, since $x\in E_{\theta}$) and $k\in[|\theta|]$. The push forward of the Lebesgue measure $\lambda_{|\theta|}$ on $E^{|\theta|}$ under the map $R_{\theta}^{-1}$ is denoted by $\lambda_{\theta}$. We note that $\lambda_{\theta}$ and $\lambda_{\theta'}$ are singular if $\theta\neq\theta'$. Let $A_{\theta}$ be the $n\times n$-matrix defined by
$$
A_{\theta}=\diag\{A_{\theta_1},\ldots,A_{\theta_p}\},
$$
where
$$
A_{\theta_k}=\frac{1}{m_{\theta_k}}\left(\begin{array}{ccc}
                    \sqrt{m_{i_k}}&\ldots&\sqrt{m_{j_k}}\\
                    \ldots&\ldots&\ldots\\
                    \sqrt{m_{i_k}}&\ldots&\sqrt{m_{j_k}}
                   \end{array}
\right)
$$
for $\theta_k=\{i_k,\ldots,j_k\}$, $i_k<\ldots< j_k$, and $m_{\theta_k}=\sum_{i\in\theta_k}m_i$,\ \  $k\in[|\theta|]$. 

We say that $f:E^n\to\R$ belongs to $C_0^2(E^n)$ if it has a compact support and can be extended to a twice continuously differentiable function $\widetilde{f}$ on an open set that contains $E^n$. Set $\frac{\partial}{\partial x_i}f(x):=\frac{\partial}{\partial x_i}\widetilde{f}(x)$, $x\in E^n$, $i\in[n]$.  Let
\begin{align}
\nabla_{\theta}f(x):=\left(\frac{1}{\sqrt{m_{\theta_k}}}\frac{\partial}{\partial y_k}f(R_{\theta}^{-1}(y))\big|_{y=R_{\theta}(x)}\right)_{k\in[|\theta|]},\quad x\in E_{\theta},
\end{align}
and 
$$
\Laplas_{\theta}f(x):=\tr\left(A_{\theta}A_{\theta}^{T}\nabla^2f\right)=\sum_{k=1}^{|\theta|}\frac{1}{m_{\theta_k}}\frac{\partial^2}{\partial y_k^2}f(R_{\theta}^{-1}(y))\big|_{y=R_{\theta}(x)},\quad x\in E_{\theta},
$$
for $f\in C_0^2(E^n)$, where $A^T$ denotes the transpose matrix.

\subsection{Construction of the finite particle system via Dirichlet form approach}\label{Dirichlet_form_for_n}

We define the measure $\Lambda_n$ on $E^n$, that will play a role of an invariant measure for a system of particles, as follows
$$
\Lambda_n=\sum_{\theta\in\varTheta^n}c_{\theta}\lambda_{\theta},
$$
where $c_{\theta}$, $\theta\in\varTheta^n$, are positive constants that will be chosen later. We also consider the following symmetric bilinear form on $L_2(E^n,\Lambda_n)$ defined on all functions $f,g$ from $C_0^2(E^n)$ by
\begin{align}
\e_n(f,g)&=\frac{1}{2}\sum_{\theta\in\varTheta^n}\int_{E^n}\langle\nabla_{\theta}f(x),\nabla_{\theta}g(x)\rangle_{\R^{|\theta|}}\Lambda_n(dx)\\
&=\frac{1}{2}\sum_{\theta\in\varTheta^n}c_{\theta}\int_{E^{|\theta|}}\left(\sum_{k=1}^{|\theta|}\frac{\partial}{\partial y_k}f(R_{\theta}^{-1}(y))\frac{\partial}{\partial y_k}g(R_{\theta}^{-1}(y))\frac{1}{m_{\theta_k}}\right)\lambda_{|\theta|}(dy),
\end{align}
where $\langle x,y\rangle_{\R^p}=\sum_{k=1}^px_ky_k$.

For each $\theta\in\varTheta^n$ we denote  
$$
\partial\theta=\left\{\theta'\in\varTheta^n:\
\begin{array}{r}
     \theta'=(\theta_1,\ldots,\theta_{k-1},\theta_k\cup\theta_{k+1},\theta_{k+2},\ldots,\theta_{|\theta|}) \\
     \mbox{for some}\ \ k\in[|\theta|-1] 
\end{array}
\right\}
$$
and define for $\theta'=(\theta_j')\in\partial\theta$ the vector $b^{\theta,\theta'}\in\R^n$ as follows
$$
b^{\theta,\theta'}_i=\begin{cases}
     -\frac{1}{m_{\theta_k}},& i\in\theta_k,\\
     \frac{1}{m_{\theta_{k+1}}},& i\in\theta_{k+1},\\
     0,& \mbox{otherwise},\\
    \end{cases}\quad i\in[n],
$$
where $k$ satisfies $\theta_k\cup\theta_{k+1}=\theta_k'$.

Using integration by parts formula, it is easily to prove the following statement.

\begin{lemma}\label{lemma_integration_by_parts_finite_dim_case}
 For each $f,g\in C_0^2(E^n)$ the relation
 $$
 \e_n(f,g)=-\int_{E^n}L_nf(x)g(x)\Lambda_n(dx)
 $$
 holds, where 
 $$
 L_nf(x)=\frac{1}{2}\sum_{\theta\in\varTheta^n}\Laplas_{\theta}f(x)\I_{E_{\theta}}(x)+\frac{1}{2}\sum_{\theta\in\varTheta^n}\langle b^{\theta},\nabla f(x)\rangle\I_{E_{\theta}}(x)
 $$
 and
 $$
 b^{\theta}=\frac{1}{c_{\theta}}\sum_{\tilde{\theta}: \theta\in\partial\tilde{\theta}}c_{\tilde{\theta}}b^{\tilde{\theta},\theta}.
 $$
\end{lemma}

It is obvious that $(L_n,C_0^2(E^n))$ is a non negative symmetric linear operator on $L_2(E^n,\Lambda_n)$. Consequently, the bilinear form $(\e_n,C_0^2(E^n))$ is closable, by~\cite[Proposition~I.3.3]{Ma:1992}. We will denote its closure by $(\e_n,\Dom_n)$.

\begin{theorem}
 $(i)$ The bilinear form $(\e_n,\Dom_n)$ is a densely defined, local, regular, conservative, symmetric Dirichlet form on $L_2(E^n,\Lambda_n)$.
 
 \vspace{7pt}
 $(ii)$ There exists a conservative diffusion process\footnote{see~\cite[Definition~V.1.10]{Ma:1992}}, i.e. a strong Markov process with continuous sample paths and infinite life time,  
 $$
 X^n=\left(\Omega^n,\F^n,(\F_t^n)_{t\geq 0},(X_t^n)_{t\geq 0}, (\p_x^n)_{x\in E^n}\right)
 $$ 
 with state space $E^n$ that is properly associated with $(\e_n,\Dom_n)$.
 
 \vspace{7pt}
 $(iii)$  The process $X^n$ is a weak solution to the SDE
 \begin{align}\label{f_SDE}
  \begin{split}
  dX^n_t&=\sum_{\theta\in\varTheta^n}A_{\theta}\I_{E_{\theta}}(X^n_t)dw(t)+\frac{1}{2}\sum_{\theta\in\varTheta^n}b^{\theta}\I_{E_{\theta}}(X^n_t)dt,\\
  X^n_0&=x
  \end{split}
 \end{align}
  in $E^n$ for $\e_n$-q.e. $x\in E^n$, where $w(t)$, $t\geq 0$, is an $n$-dimensional standard Brownian motion. 
\end{theorem}

\begin{proof}
 The proof of theorem follows from the standard arguments (see e.g.~\cite[Section~3]{Grothaus:2014}). 
\end{proof}

Choosing constants $c_{\theta}$, $\theta\in\varTheta^n$, by a special way, we can simplify equation~\eqref{f_SDE}. Let $P_{\theta}$ be the matrix defined similarly as $A_{\theta}$ with $\sqrt{m_i}$ replaced by $m_i$ for all $i\in[n]$. 

\begin{remark}
If the space $\R^n$ is furnished with the inner product $\langle x,y\rangle=\sum_{i=1}^nx_iy_im_i$, $x,y\in\R^n$, then the linear operator
$$
x\to P_{\theta}x,\quad x\in \R^n,
$$
is the orthogonal projection on $\R_{\theta}$, where $\R_{\theta}\subseteq\R^n$ is defined similarly as $E_{\theta}$ with $E^n$ replaced by $\R^n$.
\end{remark}

We also set $P_x:=P_{\theta}$ for each $x\in E_{\theta}$.

\begin{proposition}\label{prop_proces_x_n}
 Let $\varsigma\in E^n_0$. If 
 \begin{equation}\label{f_c_theta}
 c_{\theta}=\left(\prod_{k=1}^{|\theta|}m_{\theta_k}\right)\left(\prod_{k=1}^{|\theta|-1}(\varsigma_{i_k^{\theta}+1}-\varsigma_{i_k^{\theta}})\right),\quad \theta\in\varTheta^n,
 \end{equation}
 where $i_k^{\theta}=\max\theta_k$, then $b^{\theta}=\varsigma-P_{\theta}\varsigma$. Moreover, the process $X$ is a weak solution in $E^n$ to the stochastic differential equation
 \begin{align}\label{f_SDE1}
  \begin{split}
  dX^n_t&=P_{X^n_t}dB(t)+\frac{1}{2}(\varsigma-P_{X^n_t}\varsigma)dt,\\
  X^n_0&=x
  \end{split}
 \end{align}
 for $\e_n$-q.e. $x\in E^n$, where $B(t)$, $t\geq 0$, is an $n$-dimensional Brownian motion with 
$$
\Var(B_i(t),B_j(t))=\frac{t}{m_i}\I_{\{i=j\}},\quad i,j\in[n].
$$ 
\end{proposition}

\begin{proof}
We first show that $b^{\theta}=\varsigma-P_{\theta}\varsigma$. Let $\theta\in\varTheta^n$ be fixed. We will suppose that $\theta\neq(\{i\})_{i\in[n]}$, since the case $\theta=(\{i\})_{i\in[n]}$ is trivial. We also fix $j\in[n]$ and take $k$ such that $j\in\theta_k$.

Let
$$
\underline{j}:=\min\theta_k,\quad\overline{j}:=\max\theta_k
$$
and for each $l\in\{\underline{j},\ldots,\overline{j}-1\}$ we denote the sets $\{\underline{j},\ldots,l\}$ and $\{l+1,\ldots,\overline{j}\}$ by $\{\leq l\}$ and $\{>l\}$, respectively. Since $b_j^{\tilde{\theta},\theta}=0$ for all $\tilde{\theta}\in\varTheta^n$ satisfying $\theta\in\partial\tilde{\theta}$ and $\tilde{\theta}_k\cup\tilde{\theta}_{k+1}\neq\theta_k$, it is easy to see that
$$
b_j^{\theta}=\begin{cases}
               \frac{1}{c_{\theta}}\sum_{l=\underline{j}}^{\overline{j}-1}c_{\theta^l}b^{\theta^l,\theta}_j,&\underline{j}<\overline{j},\\
               0,&\underline{j}=\overline{j},
             \end{cases}
$$
where $\theta\in\partial\theta^l$ with $\theta^l_k=\{\leq l\}$ and $\theta^l_{k+1}=\{> l\}$.
We assume that $\underline{j}<\overline{j}$, otherwise $b^{\theta}_j=\varsigma_j-(P_{\theta}\varsigma)_j=0$. The simple computation gives
$$
\frac{c_{\theta^l}}{c_{\theta}}=\frac{m_{\{\leq l\}}m_{\{>l\}}}{m_{\theta_k}}(\varsigma_{l+1}-\varsigma_l)
$$ 
and
$$
b_j^{\theta^l,\theta}=\begin{cases}
                      -\frac{1}{m_{\{\leq l\}}},&l\geq j,\\
                      \frac{1}{m_{\{>l\}}},&l<j,
                      \end{cases}
$$
for all $l\in\{\underline{j},\ldots,\overline{j}-1\}$. Hence,
\begin{align}
b_j^{\theta}&=\frac{1}{m_{\theta_k}}\left[\sum_{l=\underline{j}}^{j-1}m_{\{\leq l\}}(\varsigma_{l+1}-\varsigma_l)-\sum_{l=j}^{\overline{j}-1}m_{\{>l\}}(\varsigma_{l+1}-\varsigma_l)\right]\\
&=\frac{1}{m_{\theta_k}}\left[m_{\{\leq j-1\}}\varsigma_j-\sum_{l=\underline{j}}^{j-1}m_l\varsigma_l+m_{\{>j-1\}}\varsigma_j-\sum_{l=j}^{\overline{j}}m_l\varsigma_l\right]\\
&=\varsigma_j-\frac{1}{m_{\theta_k}}\sum_{l=\underline{j}}^{\overline{j}}m_l\varsigma_l=\varsigma_j-(P_{\theta}\varsigma)_j.
\end{align}
Thus, $b^{\theta}=\varsigma-P_{\theta}\varsigma$.

The equality of the diffusion parts of~\eqref{f_SDE} and~\eqref{f_SDE1} is trivial for $B_i(t)=\frac{w_i(t)}{\sqrt{m_i}}$, $i\in[n]$. The proposition is proved. 
\end{proof}

The following example shows that one cannot expect that equation~\eqref{f_SDE1} has a strong solution.

\begin{example}
Let $n=2$, $m_1=m_2=\frac{1}{2}$ and $\varsigma=(0,1)$. Then $X_t=(x_1(t),x_2(t))$, $t\geq 0$, solves the equation
\begin{align}
dx_1(t)&=\sqrt{2}\I_{\{x_1(t)\neq x_2(t)\}}dw_1(t)\\
&+\I_{\{x_1(t)=x_2(t)\}}\frac{dw_1(t)+dw_2(t)}{\sqrt{2}}-\frac{1}{4}\I_{\{x_1(t)=x_2(t)\}}dt,\\
dx_2(t)&=\sqrt{2}\I_{\{x_1(t)\neq x_2(t)\}}dw_2(t)\\
&+\I_{\{x_1(t)=x_2(t)\}}\frac{dw_1(t)+dw_2(t)}{\sqrt{2}}+\frac{1}{4}\I_{\{x_1(t)=x_2(t)\}}dt,
\end{align}
where $(w_1,w_2)$ is a 2-dimensional standard Brownian motion. Taking
$$
y_1(t)=\frac{x_2(t)-x_1(t)}{2}\quad\mbox{and}\quad y_2(t)=\frac{x_2(t)+x_1(t)}{2},\quad t\geq 0,
$$
it is easily seen that $y_1$ and $y_2$ are weak solutions to the equations
\begin{align}
dy_1(t)&=\I_{\{y_1(t)>0\}}d\tilde{w}_1(t)+\frac{1}{4}\I_{\{y_1(t)=0\}}dt,\\
dy_2(t)&=d\tilde{w}_2(t).
\end{align}
But the equation for $y_1$ has no strong solution, according to~\cite{Engelbert:2014}. 
\end{example}

\section{\texorpdfstring{$\sigma$}{S}-finite measure on \texorpdfstring{$L_2^{\uparrow}$}{L2u}}\label{section_measure}

In Proposition~\ref{prop_proces_x_n}, we constructed the conservative diffusion process $X^n_t=(X^n_{i,t})_{i\in[n]}$ which describe the evolution of a finite sticky-reflected particle system and whose invariant measure is $\Lambda_n$ with $c_\theta$ defined by~\eqref{f_c_theta}. Moreover, it is a solution to SDE~\eqref{f_SDE1}. The goal of this section is to build the measure which will play a role of invariant measure for the infinite particle system. Since the particles keep their order, we will work with the state space $L_2^{\uparrow}$ instead of $E^n$. In this case, we can identify $X^n$ with the continuous process
$$
\sum_{i=1}^nX^n_{i,t}\I_{[a_{i-1},a_i)},\quad t\geq 0,
$$
on $L_2^{\uparrow}$, where $a_i=a_{i-1}+m_i$, $i\in[n]$, and $a_0=0$. Abusing the notation, we will also denote this process by $X_t^n$. It is easy to see, that the process $X^n$ is a solution to SDE~\eqref{f_SDE_for_X_intro} with the interaction potential $\xi$ equals $\sum_{i=1}^n\varsigma_i\I_{[a_{i-1},a_i)}$. The goal of this section is to define the invariant measure $\Xi$ for the particle system in the case of an arbitrary bounded interaction potential $\xi\in D^{\uparrow}$, which would also coincide with the push forward of the measure $\Lambda_n$ under the map 
$$
E^n\ni x\mapsto \sum_{i=1}^nx_i\I_{[a_{i-1},a_i)}\in L_2^{\uparrow}
$$
for $\xi=\sum_{i=1}^n\varsigma_i\I_{[a_{i-1},a_i)}$.

Hereinafter $\xi\in D^{\uparrow}$ is a fixed bounded function.

\subsection{Motivation of the definition}

Here we will make some manipulations with the measure $\Lambda_n$ in order to guess the formula for the measure $\Xi$. Let $n\in\N$, $m_i=\frac{i}{n}$, $i\in[n]$, and the constants $c_{\theta}$ from the definition of $\Lambda_n$ be defined by~\eqref{f_c_theta} for some $\varsigma$ that will be chosen later. We find the push forward $\widetilde{\Lambda}_n$ of the measure 
$$
\Lambda_n=\sum_{\theta\in\varTheta^n}\left(\prod_{k=1}^{|\theta|}m_{\theta_k}\right)\left(\prod_{k=1}^{|\theta|-1}(\varsigma_{i_k^{\theta}+1}-\varsigma_{i_k^{\theta}})\right)\lambda_{\theta}
$$
on $E^n$ under the map
$$
x\mapsto G(x)=\sum_{i=1}^nx_i\I_{\left[\frac{i-1}{n},\frac{i}{n}\right)},\quad x\in E^n.
$$
The measure $\widetilde{\Lambda}_n$ can be written as follows
$$
\widetilde{\Lambda}_n=\sum_{\theta\in\varTheta^n}\left(\prod_{k=1}^{|\theta|}m_{\theta_k}\right)\left(\prod_{k=1}^{|\theta|-1}(\varsigma_{i_k^{\theta}+1}-\varsigma_{i_k^{\theta}})\right)\widetilde{\lambda}(m_{\theta_1},\ldots,m_{\theta_{|\theta|}}),
$$
where $\widetilde{\lambda}(m_{\theta_1},\ldots,m_{\theta_{|\theta|}})$ is the push forward of the Lebesgue measure $\lambda_{|\theta|}$ on $E^{|\theta|}$ under the map $x\mapsto\sum_{k=1}^{|\theta|}x_k\I_{\left[a_{k+1},a_k\right)}$, with $a_0=0$, $a_k=m_{\theta_k}+a_{k-1}$, $k\in[|\theta|]$.

Setting $\varTheta_p^n=\{\theta\in\varTheta^n:\ |\theta|=p\}$ and $\varsigma_{i+1}-\varsigma_i\approx\frac{1}{n}\xi'\left(\frac{i}{n}\right)$ (if $\xi$ is continuously differentiable), it is easy to see that 
\begin{align}
 \widetilde{\Lambda}_n=&\sum_{p=1}^n\sum_{\theta\in\varTheta_p^n}\left[\prod_{k=1}^{p}\frac{|\theta_k|}{n}\right]\left[\prod_{k=1}^{p-1}\xi'\left(\frac{i_k^{\theta}}{n}\right)\frac{1}{n}\right]\widetilde{\lambda}\left(\frac{|\theta_1|}{n},\ldots,\frac{|\theta_p|}{n}\right)\\
 =&\sum_{p=1}^n\sum_{\begin{array}{l}
                      l_1,\ldots,l_p\geq 1\\
                      l_1+\ldots+l_{p}=n
                     \end{array}}
\left[\prod_{k=1}^p\frac{l_k}{n}\right]\frac{1}{n^{p-1}}\left[\prod_{k=1}^{p-1}\xi'\left(\frac{l_1+\ldots+l_k}{n}\right)\right]\widetilde{\lambda}\left(\frac{l_1}{n},\ldots,\frac{l_p}{n}\right).
\end{align}

Thus, we see that the relation consist of Riemann sums. Therefore, $\widetilde{\Lambda}_n$ looks line
\begin{align}
\sum_{p=1}^n&\int_{\begin{array}{l}
                      r_1,\ldots,r_{p-1}> 0\\
                      r_1+\ldots+r_{p-1}<1
                     \end{array}}
\left(\prod_{k=1}^{p-1}r_k\right)\left(1-r_1-\ldots-r_{p-1}\right)\\
&\cdot\left(\prod_{k=1}^{p-1}\xi'\left(r_1+\ldots+r_k\right)\right)\widetilde{\lambda}\left(r_1,\ldots,r_{p-1},1-r_1-\ldots-r_{p-1}\right)dr\\
&=\sum_{p=1}^n\int_{ 0<q_1<\ldots<q_{p-1}<1}\left(\prod_{k=1}^p(q_k-q_{k-1})\right)\left(\prod_{k=1}^{p-1}\xi'(q_k)\right)\widetilde{\lambda}\left(q_1,q_2-q_1\ldots,1-q_{p-1}\right)dq\\
&=\sum_{p=1}^n\int_{ 0<q_1<\ldots<q_{p-1}<1}\left(\prod_{k=1}^p(q_k-q_{k-1})\right)\widetilde{\lambda}\left(q_1,q_2-q_1\ldots,1-q_{p-1}\right)d\xi^{\otimes (p-1)}(q),
\end{align}
where $q_0=0$ and $q_p=1$ in the product.

In the next section, we will use the obtained expression for the definition of the measure $\Xi$. 

\subsection{Definition of the invariant measure on \texorpdfstring{$\Li$}{Li}}

We first define a measure $\Xi_n$ on $L_2$ for each $n\in\N$, supported on step functions with at most $n-1$ jumps. Let $\chi_n:Q^n\times E^n\to\Li$ with
\begin{equation}\label{f_def_of_chi}
\chi_n(q,x)=\sum_{i=1}^nx_i\I_{[q_{i-1},q_i)}+x_n\I_{\{1\}},\quad x\in E^n,\ \ q\in Q^n,
\end{equation}
and
$$
\chi_1(x)=x\I_{[0,1]},\quad x\in\R.
$$
Denote for all $q\in Q^n$, $n\geq 2$, the push forward of the Lebesgue measure $\lambda_n$ on $E^n$ under the map $\chi_n(q,\cdot)$ by $\nu_n(q,\cdot)$, i.e.
$$
\nu_n(q,A)=\lambda_n\{x:\ \chi_n(q,x)\in A\},\quad A\in\B(\Li),
$$
and set 
$$
\Xi_n(A)=\int_{Q^n}\left(\prod_{i=1}^n(q_i-q_{i-1})\right)\nu_n(q,A)d\xi^{\otimes (n-1)}(q),\ \ A\in\B(\Li),
$$
where $\int_{Q^n}\ldots d\xi^{\otimes (n-1)}(q)$ is the $(n-1)$-dim Lebesgue-Stieltjes integral with respect to $\xi^{\otimes (n-1)}(q)=\xi(q_1)\cdot\ldots\cdot\xi(q_{n-1})$. We also set
\begin{equation}\label{f_Xi_1}
\Xi_1(A)=\lambda_1\left\{x\in\R:\ \chi_1(x)\in A\right\},\quad A\in\B(\Li).
\end{equation}

Now we define the measure on $\Li$, that will be used for the definition of the Dirichlet form, as a sum of $\Xi_n$, that is,
\begin{equation}\label{f_definition_of_measure_Xi}
\Xi:=\sum_{n=1}^{\infty}\Xi_n=\Xi_1+\sum_{n=2}^{\infty}\int_{Q^n}\left(\prod_{i=1}^n(q_i-q_{i-1})\right)\nu_n(q,\cdot)d\xi^{\otimes (n-1)}(q).
\end{equation}

\begin{remark}
 If $\xi=\chi_n(q,\varsigma)$ for some $q\in Q^n$ and $\varsigma\in E^n_0$, then a simple calculation shows that $\Xi$ coincides with the push forward of the measure $\Lambda_n$ on $E^n$, defined in Section~\ref{Dirichlet_form_for_n}, under the map $x\mapsto\chi_n(q,x)$, $x\in E^n$, for $m_i=q_i-q_{i-1}$, $i\in[n]$, and $c_{\theta}$, $\theta\in\varTheta^n$, given by~\eqref{f_c_theta}.
\end{remark}

\subsection{Some properties of the measure \texorpdfstring{$\Xi$}{Xi}}

In this section, we prove some properties of the measures $\Xi$ and $\Xi_n$, $n\geq 1$. Define on $Q^n$ the measure $\mu_{\xi}^n$ as follows
$$
\mu_{\xi}^n(A)=\int_A\prod_{i=1}^n(q_i-q_{i-1})d\xi^{\otimes (n-1)}(q),\quad A\in\B(Q^n),\ \ n\geq 2.
$$

\begin{lemma}\label{lem_prop_of_Xi_n}
For each $n\in\N$, the following statements hold.
 \begin{enumerate}
  \item[(i)] The measure $\Xi_n$ is the push forward of the measure $\mu_{\xi}^n\otimes\lambda_n$ under the map $\chi_n$, if $n\geq 2$.

  \item[(ii)] The measure $\Xi_n$ is $\sigma$-finite on $\Li$ and
  $$
  \Xi_n(B_r)\leq\frac{2\pi^{\frac{n}{2}}r^n}{n!\Gamma\left(\frac{n}{2}\right)}(\xi(1)-\xi(0))^{n-1},
  $$
  where $B_r=\{g\in\Li:\ \|g\|_2\leq r\}$.

  \item[(iii)] $\Xi_n(\{g\in\Li:\ \sharp g\neq n\})=0$, where $\sharp g$ denotes the number of distinct values of the \cdl version of $g$.
 \end{enumerate}
\end{lemma}

\begin{remark}
 We note that $\{g\in\Li:\ \sharp g\neq n\}\in\B(\Li)$, since $\{g\in\Li:\ \sharp g\leq n\}$ is closed in $\Li$.
\end{remark}

\begin{remark}\label{rem_Xi_finite_on_ball}
 Property $(ii)$ of Lemma~\ref{lem_prop_of_Xi_n} immediately implies that $\Xi$ is a $\sigma$-finite measure on $\Li$ with $\Xi(B_r)<\infty$.
\end{remark}

\begin{proof}[Proof of Lemma~\ref{lem_prop_of_Xi_n}]
$(i)$ follows from the definition of the measure $\Xi_n$ and Fubini's theorem.

The equality $\nu_n(q,\{g\in\Li:\ \sharp g\neq n\})=0$, for all $q\in Q^n$, implies $(iii)$.

We note that $(ii)$ trivially holds for $n=1$ and prove $(ii)$ for $n\geq 2$. Let $q\in Q^n$ be fixed. We first estimate
\begin{align}\label{f_est_of_nu_B_r}
\begin{split}
 \nu_n(q,B_r)&=\lambda_n\left\{x\in E^n:\ \|\chi_n(q,x)\|_2^2\leq r^2\right\}\\
 &=\lambda_n\left\{x\in E^n:\ \sum_{i=1}^nx_i^2(q_i-q_{i-1})\leq r^2\right\}\\
 &\leq\frac{2\pi^{\frac{n}{2}}r^n}{n\Gamma\left(\frac{n}{2}\right)}\frac{1}{\sqrt{\prod_{i=1}^n(q_i-q_{i-1})}}.
\end{split}
\end{align}
Here $\lambda_n\left\{x\in E^n:\ \sum_{i=1}^nx_i^2(q_i-q_{i-1})\leq r^2\right\}$ is estimated by the $n$-dimensional volume of the ellipsoid $\sum_{i=1}^nx_i^2(q_i-q_{i-1})\leq r^2$. Thus,
\begin{align}
 \Xi_n(B_r)&\leq\frac{2\pi^{\frac{n}{2}}r^n}{n\Gamma\left(\frac{n}{2}\right)}\int_{Q^n}\sqrt{\prod_{i=1}^n(q_i-q_{i-1})}d\xi^{\otimes (n-1)}(q)\\
 &\leq\frac{2\pi^{\frac{n}{2}}r^n}{n\Gamma\left(\frac{n}{2}\right)}\int_{Q^n}1d\xi^{\otimes (n-1)}(q)=\frac{2\pi^{\frac{n}{2}}r^n}{n!\Gamma\left(\frac{n}{2}\right)}(\xi(1)-\xi(0))^{n-1}.
\end{align}
This completes the proof of the lemma. 
\end{proof}

The following lemma is important for the proof of the quasi-regularity of the Dirichlet form in Section~\ref{subsection_funct_with_compact_supp}. 

\begin{lemma}\label{lemma_estimation_of_int_of_norm_with_peoj}
 Let $C>0$, $q\in[1,2]$, $p,r\in[2,\infty)$ and $l\in[1,\infty)$ such that $\frac{l}{r}+\frac{2}{q}-\frac{l}{p}\leq\frac{3}{2}$ and $r\leq p$. Then there exists a constant $\tilde{C}$ which only depends on $C$ and $l$ such that 
 $$
 \sup_{h\in H}\int_{\Li}\|g\|_p^l\|\pr_gh\|_2^2\,\I_{\{\|g\|_r\leq C\}}\,\Xi(dg)\leq\tilde{C},
 $$
 where $H=\{h\in L_2:\ \|h\|_q\leq 1\}$.
\end{lemma}

\begin{proof}
 We first estimate $\int_{B_C}\|g\|_p^l\|\pr_gh\|_2^2\,\Xi_n(dg)$ for each $n\geq 2$ and $\|h\|_q\leq 1$, where $B_C=\{g\in\Li:\ \|g\|_r\leq C\}$. 
 
 By the definition of $\Xi_n$, we have
 \begin{align}
  \int_{B_C}&\|g\|_p^l\|\pr_gh\|_2^2\,\Xi_n(dg)=\int_{Q^n}\prod_{i=1}^n(q_i-q_{i-1})\\
  &\cdot\left[\int_{E^n}\left(\sum_{i=1}^n|x_i|^p(q_i-q_{i-1})\right)^{\frac{l}{p}}\left\|\pr_{\chi_n(q,x)}h\right\|_2^2\I_{B_C}(\chi_n(q,x))\lambda_n(dx)\right]d\xi^{\otimes (n-1)}(q).
 \end{align}
 Next, let $(q,x)\in Q^n\times E^n$ and $\chi_n(q,x)\in B_C$. Then
 $$
 \|\chi_n(q,x)\|_r^r=\sum_{i=1}^n|x_i|^r(q_i-q_{i-1})\leq C^r.
 $$
 Thus, $|x_i|\leq\frac{C}{(q_i-q_{i-1})^{\frac{1}{r}}}$, $i\in[n]$, and, consequently, 
 \begin{equation}\label{f_estim_of_norm_p}
 \|\chi_n(q,x)\|_p^p=\sum_{i=1}^n|x_i|^p(q_i-q_{i-1})\leq C^p\sum_{i=1}^n(q_i-q_{i-1})^{1-\frac{p}{r}}.
 \end{equation}

 Similarly, if $\|\pr_{\chi_n(q,x)}h\|_q\leq 1$, then
 \begin{equation}\label{f_estim_norm_of_pr}
 \|\pr_{\chi_n(q,x)}h\|_2^2\leq \sum_{i=1}^n(q_i-q_{i-1})^{1-\frac{2}{q}}.
 \end{equation}
 We note that, by~Remark~\ref{rem_proj_and_ecpect}~(i) and Jensen's inequality, we have that $\|h\|_q\leq 1$ implies $\|\pr_gh\|_q\leq 1$. Indeed,
 $$
 \|\pr_gh\|_q^q=\E\left|\E(h|\sigma^{\star}(g))\right|^q\leq\E\E(|h|^q|\sigma^{\star}(g))=\E|h|^q=\|h\|_q^q\leq 1.
 $$
 Thus,~\eqref{f_estim_norm_of_pr} holds for any $h\in H$. Hence, using the fact that $q_i-q_{i-1}\leq 1$, $i\in[n]$, and the inequalities~\eqref{f_estim_norm_of_pr},~\eqref{f_estim_of_norm_p}, we can estimate for each $h\in H$
 \begin{align}
 \prod_{i=1}^n(q_i-q_{i-1})&\left(\sum_{i=1}^n|x_i|^p(q_i-q_{i-1})\right)^{\frac{l}{p}}\left\|\pr_{\chi_n(q,x)}h\right\|_2^2\I_{B_C}(\chi_n(q,x))\\
 &\leq C^l\prod_{i=1}^n(q_i-q_{i-1})\left(\sum_{i=1}^n(q_i-q_{i-1})^{1-\frac{p}{r}}\right)^{\frac{l}{p}}\left(\sum_{i=1}^n(q_i-q_{i-1})^{1-\frac{2}{q}}\right)\\
 &\leq C^ln^{\frac{l}{p}}\prod_{i=1}^n(q_i-q_{i-1})^{\frac{1}{2}}\left(\sum_{i=1}^n(q_i-q_{i-1})^{\frac{3}{2}-\frac{l}{r}-\frac{2}{q}+\frac{l}{p}}\right)\\
 &\leq C^ln^{\frac{l}{p}+1}\prod_{i=1}^n(q_i-q_{i-1})^{\frac{1}{2}}\I_{B_C}(\chi_n(q,x)),
 \end{align}
 if $\frac{l}{r}+\frac{2}{q}-\frac{l}{p}\leq\frac{3}{2}$ and $r\leq p$.
 Hence, by~\eqref{f_est_of_nu_B_r} and the inclusion $B_C\subseteq\{g\in\Li:\ \|g\|_2\leq C\}$, $r\geq 2$, we have
 \begin{align}
  \int_{B_C}&\|g\|_p^l\|\pr_gh\|_2^2\,\Xi_n(dg)\\
  &\leq C^ln^{\frac{l}{p}+1}\int_{Q^n}\prod_{i=1}^n(q_i-q_{i-1})^{\frac{1}{2}}\left[\int_{E^n}\I_{B_C}(\chi_n(q,x))\lambda_n(dx)\right]d\xi^{\otimes (n-1)}(q)\\
  &=C^ln^{\frac{l}{p}+1}\int_{Q^n}\prod_{i=1}^n(q_i-q_{i-1})^{\frac{1}{2}}\nu_n(q,B_C)d\xi^{\otimes (n-1)}(q)\\
  &\leq \frac{2\pi^{\frac{n}{2}}C^{(n+l)}n^{\frac{l}{p}+1}}{n!\Gamma\left(\frac{n}{2}\right)}(\xi(1)-\xi(0))^{n-1}.
 \end{align}

 We note that $\sum_{n=2}^{\infty}\frac{2\pi^{\frac{n}{2}}C^{(n+l)}n^{\frac{l}{p}+1}}{n!\Gamma\left(\frac{n}{2}\right)}(\xi(1)-\xi(0))^{n-1}<\infty$ and 
 $$
 \sup_{h\in H}\int_{B_C}\|g\|_p^l\|\pr_gh\|_2^2\,\Xi_1(dg)\leq\int_{-C}^C|x|^ldx,
 $$
 since $\|g\|_p=\|g\|_2$ and $\|\pr_gh\|_2=\|\pr_gh\|_q\leq\|h\|_q\leq 1$ $\Xi_1$-a.e. Therefore, the integral $\int_{B_C}\|g\|_p^l\|\pr_gh\|_2^2\,\Xi(dg)$ is uniformly bounded on $H$ by a constant that only depends on $l$ and $C$. The lemma is proved.  
\end{proof}

\begin{lemma}\label{lem_p_norm conv_to_norm}
The equality  $\Xi\left\{g\in\Li:\ \|g\|_p^p\not\to\|g\|_2^2\ \ \mbox{as}\ \ p\downarrow 2\right\}=0$ holds.
\end{lemma}

\begin{proof}
 The proof follows from the definition of the measure $\Xi$ and the fact that for all $n\geq 2$ and $q\in Q^n$,
 \begin{align}
 \nu_n&\left(q,\left\{\chi_n(q,x):\ x\in E^n\ \ \mbox{and}\ \ \|\chi_n(q,x)\|_p^p\not\to\|\chi_n(q,x)\|_2^2,\ p\downarrow 2\right\}\right)\\
 &\qquad=\lambda_n\left\{x\in E^n:\ \sum_{i=1}^nx_i^p(q_i-q_{i-1})\not\to\sum_{i=1}^nx_i^2(q_i-q_{i-1}),\ p\downarrow 2\right\}=0.
 \end{align} 
\end{proof}

\subsection{Support of the measure \texorpdfstring{$\Xi$}{Xi}}

Recall that $\Li(\xi)$ denotes the subset of all $\sigma^{\star}(\xi)$-measurable functions from $\Li$. Let  $\mu_{\xi}$ denote the Lebesgue-Stieltjes measure on $[0,1]$ generated by the function $\xi$, that is, $\mu_{\xi}((a,b])=\xi(b)-\xi(a)$ for all $a<b$ from $[0,1]$.

\begin{proposition}\label{prop_support_of_Xi}
The support of $\Xi$ coincides with $\Li(\xi)$.
\end{proposition} 

\begin{remark}\label{rem_supp_of_Xi_for_incr_xi}
 If $\xi$ is a strictly increasing function, then $\Li(\xi)=\Li$ and, consequently, $\supp\Xi=\Li$. 
\end{remark}

To prove Proposition~\ref{prop_support_of_Xi}, we will need several auxiliary lemmas.

\begin{lemma}\label{lem_jump_point_and_support}
 If $h\in\St\cap\Li(\xi)$ and $s$ is a jump point of $h$, then $s\in\supp\mu_{\xi}$.
\end{lemma}

\begin{proof}
Suppose that $s\notin\supp\mu_{\xi}$. Then there exists $\eps>0$ such that $\mu_{\xi}((s-\eps,s+\eps])=0$. Thus, $\xi(s-\eps)=\xi(s+\eps)$. By Proposition~\ref{prop_xi_measurability}, we have that $h(s-\eps)=h((s+\eps)-)$. But this contradicts the assumption that $s$ is a jump point of the non decreasing function $h$. 
\end{proof}

\begin{lemma}\label{lem_proj_of_step_function}
 Let $g,h\in\Li$ and $h$ is a step function. Then $\pr_gh$ is also a step function.
\end{lemma}

\begin{proof}
 The proof if given in the appendix. 
\end{proof}

\begin{proof}[Proof of Proposition~\ref{prop_support_of_Xi}]
{\it Step I.} First we show that $\Li(\xi)\subseteq\supp\Xi$. 

Let $g\in\Li(\xi)$ and $\eps>0$. We need to show that $\Xi(B_{\eps}(g))>0$, where $B_{\eps}(g)=\{h\in \Li:\ \|g-h\|_2<\eps\}$. Since the set of all step functions $\St$ is dense in $\Li$, there exists $\widetilde{h}\in\St$ such that $\|g-\widetilde{h}\|_2<\eps$. Hence, 
\begin{equation}\label{f_g_prh}
\|g-\pr_{\xi}\widetilde{h}\|_2=\|\pr_{\xi}(g-\widetilde{h})\|_2\leq\|g-\widetilde{h}\|_2<\eps.
\end{equation}
Setting $\overline{h}=\pr_{\xi}\widetilde{h}$ and using Lemma~\ref{lem_proj_of_step_function}, we have that $\overline{h}$ is a step function that belongs to $B_{\eps}(g)\cap\Li(\xi)$. By Remark~\ref{rem_repres_of_step_functions}, there exist $n\in\N$, $r\in Q^n$ (if $n\geq 2$) and $y\in E^n_0$ such that
$$
\overline{h}=\sum_{i=1}^ny_i\I_{[r_{i-1},r_i)}+y_n\I_{\{1\}}.
$$

If $n=1$, then it is easy to see that $\Xi_1(B_{\eps}(g))>0$. This implies $\Xi(B_{\eps}(g))>0$. Thus, we will assume that $n\geq 2$.  Using the continuity of the map $F:Q^n\times E^n_0\to\R$ given by
$$
F(q,x)=\|g-\chi_n(q,x)\|_2^2=\sum_{i=1}^n\int_{q_{i-1}}^{q_i}(g(s)-x_i)^2ds,\quad (q,x)\in Q^n\times E^n_0,
$$
where $\chi_n$ is defined by~\eqref{f_def_of_chi}, and the inequality $F(r,y)<\eps^2$ which follows from~\eqref{f_g_prh}, we can conclude that there exist neighbourhoods of $r$ and $y$ defined by
$$
R=\{q\in Q^n:\ |q_i-r_i|<\delta,\ i\in[n-1]\},\quad Y=\{x\in \R^n:\ |x_i-y_i|<\delta,\ i\in[n]\}
$$
such that $Y\subset E^n$, $\prod_{i=1}^n(q_i-q_{i-1})\geq\delta$ and $F(q,x)<\eps^2$ for all $(q,x)\in R\times Y$.
Thus, trivially, $\chi_n(q,x)\in B_{\eps}(g)$ for all $(q,x)\in R\times Y$. Therefore, we can estimate $\Xi_n(B_{\eps}(g))$ from below as follows
\begin{align}
 \Xi_n(B_{\eps}(g))&=\int_{Q^n}\prod_{i=1}^n(q_i-q_{i-1})\left(\int_{E_n}\I_{\{x:\ \chi_n(q,x)\in B_{\eps}(g)\}}\lambda_n(dx)\right)d\xi^{\otimes (n-1)}(q)\\
 &\geq\delta\int_R\left(\int_{Y}1\lambda_n(dx)\right)d\xi^{\otimes (n-1)}(q)=\delta^{n+1}\prod_{i=1}^{n-1}\mu_{\xi}((r_i-\delta,r_i+\delta)).
\end{align}
Since $\overline{h}$ belongs to $\St\cap\Li(\xi)$ and $r_i$, $i\in[n-1]$, are its jump points, 
$$
\prod_{i=1}^{n-1}\mu_{\xi}((r_i-\delta,r_i+\delta))>0,
$$ 
by Lemma~\ref{lem_jump_point_and_support}. Hence $\Xi(B_{\eps}(g))>0$ and, consequently, $\Li(\xi)\subseteq\supp\Xi$.

{\it Step II.} We will establish that for all $g\in\Li\setminus\Li(\xi)$ there exists $\eps>0$ such that $\Xi(B_{\eps}(g))=0$. 

Let $g\in\Li\setminus\Li(\xi)$ be fixed. Using Proposition~\ref{prop_xi_measurability}, we can find $a,b\in[0,1]$ such that $\xi(a)=\xi(b)$ and $g(a)<g(b-)$. Thus, for some $\delta\in(0,b-a)$
$$
g(a)<g(b-\delta)\leq g(b).
$$
This inequality and the right continuity of $g$ imply that $g$ is not a constant a.e. on $[a,b]$. 

Next we claim that there exists $\eps>0$ such that
\begin{equation}\label{f_inclusion_for_B}
B_{\eps}(g)\cap\Li\subseteq\{h\in\Li:\ h(a)<h(b)\}.
\end{equation}
Indeed, if for any $\eps>0$ we can find $h$ from $B_{\eps}(g)\cap\Li$ that is a constant on $[a,b]$, then 
\begin{align}
\eps&>\|g-h\|_2^2=\int_0^1(g(s)-h(s))^2ds\geq\int_a^b(g(s)-h(a))^2ds\\
&\geq\int_a^b\left(g(s)-\frac{1}{b-a}\int_a^bg(r)dr\right)^2ds=\eps_0>0,
\end{align}
because $g$ is not a constant a.e. on $[a,b]$. But this is impossible for $\eps\leq\eps_0$. Hence,~\eqref{f_inclusion_for_B} holds for every $\eps\leq\eps_0$.

Using the inclusion~\eqref{f_inclusion_for_B}, the get
\begin{align}
 \Xi(B_{\eps}(g))&=\Xi\left(\{h\in\Li:\ h(a)<h(b)\}\cap B_{\eps}(g)\right)\\
 &=\sum_{n=2}^{\infty}\int_{Q^n}\prod_{i=1}^n(q_i-q_{i-1})\nu_n\left(q,\{h\in\Li:\ h(a)<h(b)\}\cap B_{\eps}(g)\right)d\xi^{\otimes (n-1)}(q)
\end{align}
for all $\eps \leq\eps_0$. Let $n\geq 2$ and
$$
Q_{a,b}^n:=\{q\in Q^n:\ q_i\notin(a,b]\ \ \mbox{for all}\ \ i\in[n-1]\}.
$$
Then for all $q\in Q^n_{a,b}$ 
$$
\nu_n\left(q,\{h\in\Li:\ h(a)<h(b)\}\cap B_{\eps}(g)\right)=0,
$$
since $\nu_n(q,\cdot)$ is supported on the set of step functions that have no jumps on $(a,b]$. Moreover, due to the inclusion $Q^n\setminus Q_{a,b}^n\subseteq\bigcup_{i=1}^{n-1}Q_{a,b,i}^n$, where $Q_{a,b,i}^n:=\{q\in[0,1]^{n-1}:\ q_i\in(a,b]\}$, and the equality $\xi(a)=\xi(b)$, we have
\begin{align}
\mu_{\xi}^n(Q^n\setminus Q_{a,b}^n)&\leq\sum_{i=1}^{n-1}\mu_{\xi}^n(Q_{a,b,i}^n)=
\int_{Q_{a,b,i}^n}\prod_{i=1}^n(q_i-q_{i-1})d\xi^{\otimes (n-1)}(q)\\
&\leq\sum_{i=1}^{n-1}(\xi(1)-\xi(0))^{n-2}(\xi(b)-\xi(a))=0.
\end{align}
Thus, $\Xi(B_{\eps}(g))=0$.
This finishes the proof of the proposition. 
\end{proof}

\begin{corollary}\label{corol_supp_of_Xi_n}
 If $\sharp\xi\geq n$, then $\supp\Xi_n=\Li(\xi)\cap\{g\in\Li:\ \sharp g\leq n\}$. Otherwise, $\Xi_n=0$.
\end{corollary}

\begin{proof}
 The inclusion $\supp\Xi_n\subseteq\Li(\xi)\cap\{g\in\Li:\ \sharp g\leq n\}$ immediately follows from Proposition~\ref{prop_support_of_Xi} and Lemma~\ref{lem_prop_of_Xi_n}~$(iii)$.
 
 Next assuming $\sharp\xi\geq n$, we will prove that 
 \begin{equation}\label{f_supp_of_Xi_n}
  \Xi_n(B_{\eps}(g))>0
 \end{equation}
 for all $g\in \Li(\xi)\cap\{g\in\Li:\ \sharp g\leq n\}$ and $\eps>0$. Since the close of $\{g\in\Li:\ \sharp g=n\}\cap\Li(\xi)$ in $\Li$ coincides with $\{g\in\Li:\ \sharp g\leq n\}\cap\Li(\xi)$,
 it is enough to check the inequality~\eqref{f_supp_of_Xi_n} for functions of the form
 $$
 g=\chi(q,x),\quad (q,x)\in Q^n\times E^n_0.
 $$
 Thus, fixing $g=\chi(q,x)$ for some $(q,x)\in Q^n\times E^n_0$, similarly to Step~I of the proof of Proposition~\ref{prop_support_of_Xi}, we can show that $\Xi_n(B_{\eps}(g))>0$. Hence, $\supp\Xi_n=\Li(\xi)\cap\{g\in\Li:\ \sharp g\leq n\}$.
 
 If $\sharp\xi<n$, then $\Li(\xi)\cap\{g\in\Li:\ \sharp g= n\}=\emptyset$, by Proposition~\ref{prop_xi_measurability}. Consequently, Proposition~\ref{prop_support_of_Xi} and Lemma~\ref{lem_prop_of_Xi_n}~$(iii)$ yield the equality $\Xi_n=0$. 
\end{proof}

\begin{corollary}\label{corol_full_measure_set_for_Xi}
 The set $\St\cap\Li(\xi)$ has full measure~$\Xi$, that is, $\Xi(\Li(\xi)\setminus\St)=0$.
\end{corollary}
\begin{proof}
 The corollary follows from the definition of the measure $\Xi$ and Corollary~\ref{corol_supp_of_Xi_n}. 
\end{proof}

\section{Definition of the Dirichlet form on \texorpdfstring{$\LL$}{L2(L2,xi)}}\label{section_dirichlet_form}

As before, we will assume that $\xi\in D^{\uparrow}$ is a bounded function and $\Xi$ is a measure on $\Li$ defined by~\eqref{f_definition_of_measure_Xi}. Since $\Xi$ is supported on the space $\Li(\xi)$, hereinafter we will work with spaces $\Li(\xi)$ and $L_2(\xi)$ instead of $\Li$ and $L_2$. Let $\LL$ or simpler $\Ll$ denote the space of $\Xi$-integrable functions on $\Li$ with the usual norm $\|\cdot\|_{\Ll}$ and the inner product $\langle\cdot,\cdot\rangle_{\Ll}$. The goal of this section is to construct the Dirichlet form on $\LL$ which will define an infinite sticky-reflected particle system with interaction potential $\xi$.

\subsection{A set of admissible functions on \texorpdfstring{$\Li(\xi)$}{Li(xi)}}

Let $\Cfb(\R^m)$ be the set of all infinitely differentiable (real-valued) functions on $\R^m$ with all partial derivatives bounded and $\Cfo(\R^m)$ be the set of functions from $\Cfb(\R^m)$ with compact
support. In this section, we will define the class of ``smooth'' integrable functions on $\Li(\xi)$. Since $\Li(\xi)\subseteq L_2(\xi)$, it is reasonable to consider functions of the form $u(\langle\cdot,h_1\rangle,\ldots,\langle\cdot,h_m\rangle)$, where $u\in\Cfb(\R^m)$ and $h_j\in L_2(\xi)$, $j\in[m]$. But in general, these functions are not integrable with respect to the measure $\Xi$. Therefore, we will need to cut off them by functions with bounded support in $\Li(\xi)$. Let $\FC$ denote the linear space generated by functions on $\Li(\xi)$ of the form
\begin{equation}\label{f_form_of_U}
U=u(\langle\cdot,h_1\rangle,\ldots,\langle\cdot,h_m\rangle)\varphi(\|\cdot\|_2^2)=u(\langle\cdot,\vec{h}\rangle)\varphi(\|\cdot\|_2^2),
\end{equation}
where $u\in\Cfb(\R^m),$ $\varphi\in\Cfo(\R)$ and $h_j\in L_2(\xi),$ $j\in[m]$.

\begin{remark}\label{rem_prop_of_FC}
 \begin{enumerate}
  \item[(i)] The set $\FC$ is an associative algebra, in particular, $U,V\in\FC$ implies $UV\in\FC$.
  
  \item[(ii)] Since each $U\in\FC$ has a bounded support, $\FC\subseteq\LL$, by Remark~\ref{rem_Xi_finite_on_ball}.
  
  \item[(iii)] For each $n\geq 2$ and $q\in Q^n$ the function $x\mapsto U(\chi_n(q,x))$ belongs to $\Cfo(E^n)$ and, similarly, $x\mapsto U(\chi_1(x))$ belongs to $\Cfo(\R)$.
 \end{enumerate}
\end{remark}

\begin{proposition}\label{prop_densely_of_FC}
 The set $\FC$ is dense in $\LL$.
\end{proposition}

\begin{proof}
The proof of the proposition follows from a standard approximation argument.
\end{proof}

\subsection{Differential operator and integration by parts formula}

In this section, we will define the differential operator $\D$ on $\FC$ will will be used for the definition of the Dirichlet form. 

For each $U\in\FC$ given by~\eqref{f_form_of_U} the \textit{differential operator} is defined by
\begin{equation}\label{f_derivativ_D}
\begin{split}
\D U(g):&=\pr_g\left[\nabla^{L_2}U(g)\right]\\
&=\varphi(\|g\|_2^2)\sum_{j=1}^m\partial_ju(\langle g,\vec{h}\rangle)\pr_gh_j
+u(\langle g,\vec{h}\rangle)\varphi'(\|g\|_2^2)2g, \quad g \in \Li(\xi),
\end{split}
\end{equation}
where $\nabla^{L_2}$ denotes the Fr\'{e}chet derivative on $L_2$ and $\partial_ju(y)=\frac{\partial}{\partial y_j}u(y)$, $y\in\R^m$. For every function $U$ from $\FC$, $DU$ is define by linearity.

A simple calculation gives the following statement.
\begin{lemma}\label{lemma_connection_D_with_D_n}
 For all $(q,x)\in Q^n\times E^n_0$, $n\geq 2$,
 $$
 \D U(\chi_n(q,x))=\sum_{i=1}^n\frac{\partial}{\partial x_i}U(\chi_n(q,x))\frac{\I_{[q_{i-1},q_i)}}{(q_i-q_{i-1})}
 $$
 and 
 $$
 \D U(\chi_1(x))=\frac{d}{d x}U(\chi_1(x))\I_{[0,1]}.
 $$
 In particular, for each $i\in[n]$
 $$
 \langle\D U(\chi_n(q,x)),\I_{[q_{i-1},q_i)}\rangle=\langle\nabla^{L_2} U(\chi_n(q,x)),\I_{[q_{i-1},q_i)}\rangle=\frac{\partial}{\partial x_i}U(\chi_n(q,x)).
 $$
\end{lemma}

The definition of the differential operator and Lemma~\ref{lemma_connection_D_with_D_n} imply the following trivial properties of~$\D$.

\begin{remark}\label{ref_prop_of_D}
\begin{enumerate}
 \item[(i)] For each $U\in\FC$, $\D U$ maps $\Li(\xi)$ into $L_2(\xi)$ and, in general, $\D U$ is not continuous, since $\pr_{\cdot}h$ is not, for each non constant function $h\in L_2(\xi)$.
 
 \item[(ii)] $\D$ is a linear operator satisfying the Leibniz rule.
 
 \item[(iii)] For each $U\in\FC$, $f\in L_2(\xi)$ and $g\in\Li(\xi)$
 $$
 \D_fU(g):=\langle\D U(g),f\rangle=\lim_{\eps\downarrow 0}\frac{U\left(g+\eps\pr_{g}f\right)-U(g)}{\eps}.
 $$
\end{enumerate}
\end{remark}

Now we prove the integration by parts formula for the operator $D$. For this we first define the second order differential operator on $\FC$ in a similar way as in the finite dimensional case. We set for $U\in\FC$
$$
L_0U(g)=\begin{cases}
           \sum_{i=1}^n\frac{\partial^2}{\partial x_i^2}U(\chi_n(q,x))\frac{1}{(q_i-q_{i-1})},& g=\chi_n(q,x),\ \ n\geq 2,\\
           & \quad\quad (q,x)\in Q^n\times E_0^n,\\
           \frac{d^2}{d x^2}U(\chi_1(x)),& g=\chi_1(x),\ \ x\in\R,\\
           0,&\mbox{otherwise}.
          \end{cases}
$$

Using simple computations and Remark~\ref{rem_repres_of_step_functions}, we can prove the following lemma.
\begin{lemma}
If $U\in\FC$ is given by~\eqref{f_form_of_U}, then
\begin{align}
L_0U(g)&=\varphi(\|g\|_2^2)\sum_{i,j=1}^m\partial_i\partial_ju(\langle g,\vec{h}\rangle)\langle\pr_gh_i,\pr_gh_j\rangle\\
&+u(\langle g,\vec{h}\rangle)\left[4\varphi''(\|g\|_2^2)\|g\|_2^2+2\varphi'(\|g\|_2^2)\cdot\sharp g\right]\\
&+2\sum_{j=1}^m\partial_ju(\langle g,\vec{h}\rangle)\varphi'(\|g\|_2^2)\langle\pr_gh_i,g\rangle,\quad g\in\St,
\end{align}
and 
$$
L_0U(g)=0,\quad g\in\Li(\xi)\backslash\St.
$$
\end{lemma}

\begin{theorem}[Integration by parts formula]\label{theorem_integratin_by_parts}
 Let $U,V\in\FC$. Then
 \begin{align}
 \int_{\Li(\xi)}\langle\D U(g),\D V(g)\rangle\Xi(dg)&=-\int_{\Li(\xi)}L_0U(g)V(g)\Xi(dg)\\
 &-\int_{\Li(\xi)}V(g)\langle \nabla^{L_2}U(g)-\D U(g),\xi\rangle\Xi(dg).
 \end{align}
 In particular, if $U$ is given by~\eqref{f_form_of_U}, then
 \begin{equation}\label{f_int_by_parts}
 \begin{split}
 \int_{\Li(\xi)}\langle\D U(g)&,\D V(g)\rangle\Xi(dg)=-\int_{\Li(\xi)}L_0U(g)V(g)\Xi(dg)\\
 &-\int_{\Li(\xi)}\varphi(\|g\|_2^2)V(g)\sum_{j=1}^m\partial_ju(\langle g,\vec{h}\rangle)\langle h_j,\xi-\pr_{g}\xi\rangle\Xi(dg).
 \end{split}
 \end{equation}
\end{theorem}

We remark that $\nabla^{L_2}U(g)-\D U(g)$ coincides with the projection of $\nabla^{L_2}U(g)$ onto the orthogonal complement of $L_2(g)$ in $L_2$.

\begin{proof}[Proof of Theorem~\ref{theorem_integratin_by_parts}]
 To prove the proposition, we will use Lemma~\ref{lemma_connection_D_with_D_n} and the integration by parts formula for the Riemann integral. 
 
 We first note that 
 \begin{equation}\label{f_int_by_parts_for_Xi_1}
 \int_{\Li(\xi)}\langle\D U(g),\D V(g)\rangle\Xi_1(dg)=-\int_{\Li(\xi)}L_0U(g)V(g)\Xi_1(dg).
 \end{equation}
 Indeed, by~\eqref{f_Xi_1} and Remark~\ref{rem_prop_of_FC}~(iii),
 \begin{align}
  \int_{\Li(\xi)}\langle\D U(g)&,\D V(g)\rangle\Xi_1(dg)=\int_{\R}\langle\D U(\chi_1(x)),\D V(\chi_1(x))\rangle dx\\
  &=\int_{\R}\frac{d}{dx} U(\chi_1(x)) \frac{d}{dx}V(\chi_1(x)))dx\\
  &=-\int_{\R}\left(\frac{d^2}{dx^2} U(\chi_1(x))\right)V(\chi_1(x)))dx\\
  &=-\int_{\Li(\xi)}L_0U(g)V(g)\Xi_1(dg).
 \end{align}
 
 Next, we check that for each $n\geq 2$
 \begin{align}\label{f_int_by_parts_for_Xi_n}
 \begin{split}
 \int_{\Li(\xi)}\langle\D U(g)&,\D V(g)\rangle\Xi_n(dg)=-\int_{\Li(\xi)}L_0U(g)V(g)\Xi_n(dg)\\
 &-\int_{\Li(\xi)}\langle \nabla^{L_2}U(g)-\D U(g),\xi\rangle V(g)\Xi_{n-1}(dg).
 \end{split}
 \end{align}
 To show this, we reduce the integral with respect to $\Xi_n$ to the Riemann-Stieltjes integral similarly to the previous case. Thus, by Lemma~\ref{lem_prop_of_Xi_n}~(i), we have
 \begin{align}
  &\int_{\Li(\xi)}\langle\D U(g),\D V(g)\rangle\Xi_n(dg)\\
  &=\int_{Q^n}\prod_{i=1}^n(q_i-q_{i-1})\left[\int_{E^n}\langle\D U(\chi(q,x)),\D V(\chi(q,x))\rangle\lambda_n(dx)\right]d\xi^{\otimes (n-1)}(q).
 \end{align}
 We fix $q\in Q^n$ and apply to the integral with respect to $\lambda_n$ the usual integration by parts formula. Hence, using Lemma~\ref{lemma_connection_D_with_D_n}, we obtain
 \begin{align}
  \int_{E^n}&\langle\D U(\chi(q,x)),\D V(\chi(q,x))\rangle\lambda_n(dx)\\
  &=\int_{E^n}\sum_{i=1}^n\frac{\partial}{\partial x_i}U(\chi(q,x))\frac{\partial}{\partial x_i}V(\chi(q,x))\frac{1}{q_i-q_{i-1}}\lambda_n(dx)\\
  &=-\int_{E^n}\sum_{i=1}^n\left(\frac{\partial^2}{\partial x_i^2}U(\chi(q,x))\right)\frac{1}{q_i-q_{i-1}}V(\chi(q,x))\lambda_n(dx)\\
  &+\sum_{i=1}^n\int_{E^{n-1}}\left.\left[\left(\frac{\partial}{\partial x_i}U(\chi(q,x))\right)V(\chi(q,x))\right]\right|_{x_i=x_{i-1}}^{x_i=x_{i+1}}\frac{\lambda_{n-1}(dx^{(i)})}{q_i-q_{i-1}}\\
  &=:I_1(q)+I_2(q),
 \end{align}
 where $x^{(i)}=(x_1,\ldots,x_{i-1},x_{i+1},\ldots,x_n)$, $x_0=-\infty$ and $x_{n+1}=+\infty$. 
 
 By the definition of the operator $L_0$ and Lemma~\ref{lem_prop_of_Xi_n}~(i), we have that 
 $$
 \int_{Q^n}\left(\prod_{i=1}^n(q_i-q_{i-1})\right)I_1(q)d\xi^{\otimes (n-1)}(q)=-\int_{\Li(\xi)}L_0U(g)V(g)\Xi_n(dg).
 $$
 
 Next we rewrite $I_2(q)$. By Lemma~\ref{lemma_connection_D_with_D_n}, we obtain
 \begin{align}
  I_2(q)&=\sum_{i=1}^n\int_{E^{n-1}}\left.\left[\left\langle\nabla^{L_2}U(\chi(q,x)),\I_{[q_{i-1},q_i)}\right\rangle V(\chi(q,x))\right]\right|_{x_i=x_{i-1}}^{x_i=x_{i+1}}\frac{\lambda_{n-1}(dx^{(i)})}{q_i-q_{i-1}}\\
&=\sum_{i=1}^{n-1}\int_{E^{n-1}}\left\langle\nabla^{L_2}U(\chi(q^{(i)},x)),e_i(q)-e_{i+1}(q)\right\rangle V(\chi(q^{(i)},x))\lambda_{n-1}(dx),
 \end{align}
where $q^{(i)}$ is defined similarly to $x^{(i)}$, removing the $i$-th coordinate, and $e_i(q):=\frac{\I_{[q_{i-1},q_i)}}{q_i-q_{i-1}}$, $i\in[n]$.
For simplicity of notation, we set 
$$
c_n(q)=\prod_{i=1}^n(q_i-q_{i-1}).
$$ 
Then
\begin{align}
 \int_{Q^n}c_n(q)I_2(q)&d\xi^{\otimes (n-1)}(q)=\sum_{i=1}^n\int_{E^{n-1}}\bigg[\int_{Q^n}c_n(q)\Big\langle\nabla^{L_2}U(\chi(q^{(i)},x)),e_i(q)\\&
 -e_{i+1}(q)\Big\rangle V(\chi(q^{(i)},x))d\xi^{\otimes (n-1)}(q)\bigg]\lambda_{n-1}(dx)\\
 &=\sum_{i=1}^n\int_{E^{n-1}}\bigg[\int_{Q^{n-1}}c_{n-1}(q^{(i)})\left\langle\nabla^{L_2}U(\chi(q^{(i)},x)),f(q^{(i)})\right\rangle\\
 & \cdot V(\chi(q^{(i)},x))d\xi^{\otimes (n-2)}(q^{(i)})\bigg]\lambda_{n-1}(dx),
\end{align}
where 
$$
f(q^{(i)}):=\int_{q_{i-1}}^{q_{i+1}}\frac{(q_{i+1}-q_i)(q_i-q_{i-1})}{q_{i+1}-q_{i-1}}(e_i(q)-e_{i+1}(q))d\xi(q_i).
$$

Integrating by parts, we obtain
\begin{align}
f(q^{(i)})(r)&=\left(\int_r^{q_{i+1}}\frac{q_{i+1}-q_i}{q_{i+1}-q_{i-1}}d\xi(q_i)-\int_{q_{i-1}}^r\frac{q_i-q_{i-1}}{q_{i+1}-q_{i-1}}d\xi(q_i)\right)\I_{[q_{i-1},q_{i+1})}(r)\\
&=\left(\frac{1}{q_{i+1}-q_{i-1}}\left\langle\xi,\I_{[q_{i-1},q_{i+1})}\right\rangle-\xi(r)\right)\I_{[q_{i-1},q_{i+1})}(r),\quad r\in[0,1].
\end{align}
Hence,
\begin{align}
 \int_{Q^n}c_n(q)I_2(q)d\xi^{\otimes (n-1)}(q)&=\int_{Q^{n-1}}c_{n-1}(q)\bigg[\int_{E^{n-1}}\Big\langle\nabla^{L_2}U(\chi(q,x)),\pr_{\chi(q,\widetilde{x})}\xi-\xi\Big\rangle\\
 &\cdot V(\chi(q,x))\lambda_{n-1}(dx)\bigg]d\xi^{\otimes (n-2)}(q),
\end{align}
where $\widetilde{x}$ is any point from $E^{n-1}_0$ (note that $\pr_{\chi(q,\widetilde{x})}=\pr_{\chi(q,\widetilde{y})}$ for all $\widetilde{x},\widetilde{y}\in E^{n-1}_0$). This immediately implies
\begin{align}
 \int_{Q^n}\left(\prod_{i=1}^n(q_i-q_{i-1})\right)&I_2(q)d\xi^{\otimes (n-1)}(q)\\
 &=-\int_{\Li(\xi)}\langle \nabla^{L_2}U(g),\xi-\pr_g\xi\rangle V(g)\Xi_{n-1}(dg)\\
 &=-\int_{\Li(\xi)}\langle \nabla^{L_2}U(g)-\D U(g),\xi\rangle V(g)\Xi_{n-1}(dg),
\end{align}
where we have used the equality
\begin{equation}\label{f_adjoint_og_pr}
\langle \nabla^{L_2}U(g)-\D U(g),\xi\rangle=\langle \nabla^{L_2}U(g),\xi-\pr_g\xi\rangle
\end{equation}
It proves~\eqref{f_int_by_parts_for_Xi_n}.

Now, summing~\eqref{f_int_by_parts_for_Xi_1} and~\eqref{f_int_by_parts_for_Xi_n} over $n$ and using Remark~\ref{rem_Xi_finite_on_ball}, we obtain the integration by parts formula. The expression~\eqref{f_int_by_parts} easily follows from~\eqref{f_adjoint_og_pr} and the equality $\langle g,\pr_g\xi-\xi\rangle=0$.
This completes the proof of the theorem. 
\end{proof}

The same argument as in the proof of Theorem~\ref{theorem_integratin_by_parts} yields the integration by parts formula for $\D_f=\langle\D\cdot,f\rangle$, $f\in\Li(\xi)$.

\begin{proposition}
 For each $U,V\in\FC$ and $f\in L_2$
 \begin{align}
 \int_{\Li(\xi)}\left(\D_f U(g)\right) V(g)\Xi(dg)&=-\int_{\Li(\xi)}U(g)\D_fV(g)\Xi(dg)\\
 &-\int_{\Li(\xi)}U(g)V(g)\langle f,\xi-\pr_{g}\xi\rangle\Xi(dg).
 \end{align}
\end{proposition}

\subsection{The Dirichlet form \texorpdfstring{$(\e,\Dom)$}{Li}}

We define
$$
\e(U,V)=\frac{1}{2}\int_{\Li(\xi)}\langle\D U(g),\D V(g)\rangle\Xi(dg),\quad U,V\in\FC.
$$
Then $(\e,\FC)$ is a densely defined positive definite symmetric bilinear form on $\LL$, by Proposition~\ref{prop_densely_of_FC}.
The integration by parts formula implies that there exists a negative definite symmetric linear operator $L$ on $\Ll$, given by
\begin{equation}\label{f_generator}
\begin{split}
LU(g):&=\frac{1}{2}\left[L_0U(g)+\langle \nabla^{L_2}U(g)-\D U(g),\xi\rangle\right]\\
&=\frac{1}{2}\left[L_0U(g)+\varphi(\|g\|_2^2)\sum_{j=1}^m\partial_ju(\langle g,\vec{h}\rangle)\langle\xi-\pr_g\xi,h_j\rangle\right],\quad g\in\Li(\xi),
\end{split}
\end{equation}
if $U\in\FC$ is defined by~\eqref{f_form_of_U}, such that
$$
\e(U,V)=-\langle LU,V\rangle_{\Ll}.
$$
Consequently, by~\cite[Proposition~I.3.3]{Ma:1992}, $(\e,\FC)$ is closable on $\Ll$.

\begin{definition}
The closure $(\e,\FC)$ on $\Ll$ will be denoted by $(\e,\Dom)$. 
\end{definition}

\begin{remark}
We can extend the differential operator $\D$ to $\Dom$, letting
$$
\D U:=\lim_{n\to\infty}\D U_n\quad \mbox{in}\ \ \Ll,
$$ 
if $\{U_n,\ n\geq 1\}\subset\FC$ converges to $U\in\Dom$ with respect to the norm $\e_1^{\frac{1}{2}}$, where $\e_1:=\e(\cdot,\cdot)+\langle\cdot,\cdot\rangle_{\Ll}$. Then, for all $U,V\in\Dom$
\begin{equation}\label{f_dirichlet_form}
\e(U,V)=\frac{1}{2}\int_{\Li(\xi)}\langle\D U(g),\D V(g)\rangle\Xi(dg).
\end{equation}
\end{remark}

We next check that $(\e,\Dom)$ is a Dirichlet form. For this we will need an analog of the chain rule.

\begin{lemma}\label{lemma_chain_rule}
  Let $F\in C^1(\R^k)$, $F(0)=0$  
  and $U_j\in\FC$, $j\in[k]$. Then the composition $F(U)=F(U_1,\ldots,U_k)$ belongs to $\Dom$ and 
  $$
  \D F(U)(g)=\sum_{j=1}^k\partial_jF(U(g))\D U_j(g),\quad g\in\Li(\xi).
  $$
\end{lemma}

\begin{proof}
We will prove the lemma, using the approximation of $F$ by the Bernstein polynomials and the fact that $\FC$ is an associative algebra (see Remark~\ref{rem_prop_of_FC}~(i)).

Since $U_j$, $j\in[k]$, belong to $\FC$, they are bounded, i.e. there exists a constant $M$ such that $|U_j(g)|\leq M$ for all $g\in\Li(\xi)$ and $j\in[k]$. Let $P_n^M(F;\cdot)$, $n\geq 1$, be polynomials defined by~\eqref{f_polynomial_P}. Then by~Lemma~\ref{lemma_converg_of_polynomials}, 
$$
\left|P_n^M(F;U(g))-F(U(g))\right|\leq\sup_{x\in[-M,M]^k}\left|P_n^M(F;x)-F(x)\right|\I_{\supp U}(g)\to 0,
$$
as $n\to\infty$, where $\supp U:=\bigcup_{j=1}^k\supp U_j$. Hence, by Remarks~\ref{rem_Xi_finite_on_ball},~\ref{rem_prop_of_FC}~(ii) and the dominated convergence theorem, we have that $\{P_n^M(F;U)\}_{n\geq 1}$ converges to $F(U)$ in $\Ll$.

Remark~\ref{rem_prop_of_FC}~(ii) and the fact that $P_n^M(F;0)=0$ imply that $P_n^M(F;U)\in\FC$. Moreover, the Leibniz rule for $\D$ (see Remark~\ref{ref_prop_of_D}) yields
$$
\D P_n^M(F;U)(g)=\sum_{j=1}^{k}\partial_jP_n^M(F;U(g))\D U_j(g),\quad g\in\Li(\xi).
$$
Using the estimate
\begin{align}
&\left|\partial_jP_n^M(F;U(g))\D U_j(g)-\partial_jF(U(g))\D U_j(g)\right|\\
&\leq\sup_{x\in[-M,M]^k}\left|\partial_jP_n^M(F;x)-\partial_jF(x)\right||\D U_j(g)|,
\end{align}
Lemma~\ref{lemma_converg_of_polynomials} and \ the dominated \ convergence theorem, \ we \ obtain \ that $\{\D P_n^M(F;U)\}_{n\geq 1}$ converges to $\sum_{j=1}^k\partial_jF(U)\D U_j$ in $\Ll$. It completes the proof of the lemma. 
\end{proof}

\begin{corollary}\label{corol_relaxetion_of_cond_on_u}
 For each $u\in C_b^1(\R^m)$, $h_j\in L_2(\xi),$ $j\in[m]$, and $\varphi\in\Cfo(\R)$ the function $U=u(\langle\cdot,h_1\rangle,\ldots,\langle\cdot,h_m\rangle)\varphi(\|\cdot\|_2^2),$ belongs to $\Dom$ and $\D U$ is defined by~\eqref{f_derivativ_D}.
\end{corollary}

\begin{proof}
The statement of the corollary follows from Lemma~\ref{lemma_chain_rule}.
\end{proof}

The following chain rule for $\D$ easily follows from Lemma~\ref{lemma_chain_rule} and the closability of $(\e,\Dom)$.

\begin{proposition}\label{prop_chain_rule_for_D}
 Let $F\in C^1_b(\R^k)$, $F(0)=0$  
  and $U_j\in\Dom$, $j\in[k]$. Then the function $F(U)=F(U_1,\ldots,U_k)$ belongs to $\Dom$ and 
  $$
  \D F(U)(g)=\sum_{j=1}^k\partial_jF(U(g))\D U_j(g),\quad g\in\Li(\xi).
  $$
\end{proposition}

We now are ready to prove that $(\e,\Dom)$ is a Dirichlet form on $\LL$. For $U,V:\Li(\xi)\to \R$ we set
$$
U\wedge V=\min\{U,V\}\quad\mbox{and}\quad U\vee V=\max\{U,V\}.
$$

\begin{proposition}
 The bilinear form $(\e,\Dom)$ is a symmetric Dirichlet form on $\LL$, that is, for all $U\in\Dom$ the function $(U\vee 0)\wedge 1$ belongs to $\Dom$ and
 $$
 \e((U\vee 0)\wedge 1,(U\vee 0)\wedge 1)\leq\e(U,U).
 $$
\end{proposition}

\begin{proof}
 To prove the proposition, we need to show that for each $U\in\Dom$ and $\eps>0$ there exists a function $F_{\eps}:\R\to[-\eps,1+\eps]$ such that $F_{\eps}(x)=x$ for all $x\in[0,1]$, $0\leq F_{\eps}(x_2)-F_{\eps}(x_1)\leq x_2-x_1$ if $x_1\leq x_2$, $F_{\eps}(U)\in\Dom$ and 
 $$
 \limsup_{\eps\to 0}\e(F_{\eps}(U),F_{\eps}(U))\leq\e(U,U),
 $$
 according to~\cite[Proposition~I.4.7]{Ma:1992}.
 
We take  for $\eps>0$ an arbitrary non decreasing continuously differentiable function $F_{\eps}:\R\to[-\eps,1+\eps]$ such that $|F'(x)|\leq 1$, $x\in\R$, and $F_{\eps}(x)=x$ for all $x\in[0,1]$. Then it is clear that $0\leq F_{\eps}(x_2)-F_{\eps}(x_1)\leq x_2-x_1$ if $x_1\leq x_2$. By Proposition~\ref{prop_chain_rule_for_D}, $F_{\eps}(U)\in\Dom$ and 
 $$
 \limsup_{\eps\to 0}\e(F_{\eps}(U),F_{\eps}(U))=\frac{1}{2}\limsup_{\eps\to 0}\int_{\Li(\xi)}|F_{\eps}'(U(g))|^2\|\D U(g)\|_2^2\Xi(dg)\leq\e(U,U).
 $$
 This completes the proof of the proposition.
\end{proof}

\begin{lemma}\label{lemma_estim_for_e_of_max}
 Let $U,V$ in $\Dom$. Then $U\vee V\in\Dom$ and 
 \begin{equation}\label{f_estim_for_e_of_max}
 \e(U\vee V,U\vee V)\leq\e(U,U)\vee\e(V,V).
 \end{equation}
\end{lemma}

\begin{proof}
 The fact that $U\vee V\in\Dom$ follows from~\cite[Proposition~I.4.11]{Ma:1992}. Inequality~\eqref{f_estim_for_e_of_max} can be proved similarly to~\cite[Lemma~IV.4.1]{Ma:1992}. 
\end{proof}

\begin{lemma}\label{lemma_product_belongs_to_Dom}
 Let $U,V\in\Dom$ and $|U|\vee\|\D U\|_2$ be bounded $\Xi$-a.e. Then $U\cdot V\in\Dom$ and $\D(U\cdot V)=(\D U)\cdot V+U\cdot\D V$.
\end{lemma}

\begin{proof}
 The lemma follows from~\cite[Corollary~I.4.15]{Ma:1992} and Proposition~\ref{prop_chain_rule_for_D}, using an approximation (w.r.t $\e^{\frac{1}{2}}$-norm) of $V$ by bounded functions.  
\end{proof}

\section{Quasi-regularity of the Dirichlet form \texorpdfstring{$(\e,\Dom)$}{(E,D)}}\label{section_quasi_regularity}

The goal of this section is to prove that the Dirichlet form $(\e,\Dom)$ is quasi-regular. This will imply the existence of a Markov process in $\Li(\xi)$ that is properly associated with $(\e,\Dom)$.

\subsection{Functions with compact support}\label{subsection_funct_with_compact_supp}

In this section, we will show that the domain $\Dom$ of the Dirichlet form contains a rich enough subset of functions with compact support.

\begin{lemma}\label{lemma_function_with_compact_supp}
 For every $p\in\left[2,\frac{5}{2}\right]$, $g_0\in\Li(\xi)$ and $\varphi\in\Cfo(\R)$ the function $\varphi(\|\cdot-g_0\|_p^p)$ belongs to~$\Dom$. Moreover, $\D\varphi(\|\cdot-g_0\|_2^2)(g)=2\varphi'(\|g-g_0\|_2^2)\pr_g(g-g_0)$ for all $g\in\Ll$.
\end{lemma}

\begin{proof}
 For simplicity we give the proof for $g_0=0$. 
 
 Let $\{h_n\}_{n\geq 1}\subseteq L_{\infty}$ be a dense subset in $L_q$ with $\|h_n\|_q=1$, where $\frac{1}{p}+\frac{1}{q}=1$. Then 
 $$
 \|g\|_p=\sup_{n\geq 1}|\langle g,h_n\rangle|=\sup_{n\geq 1}\left|\int_0^1g(s)h_n(s)ds\right|.
 $$
 
 Next we take functions $\psi_1,\psi_2\in\Cfo(\R)$ such that $\psi_1=1$ on $[-M-1,M+1]$, $\supp\psi_1\subseteq [-2M-2,2M+2]$, $\psi_2=1$ on $[-M,M]$ and $\supp\psi_2\subseteq [-M-1,M+1]$, where $M$ is chosen such that the interval $[-M^{\frac{p}{2}},M^{\frac{p}{2}}]$ contains $\supp\varphi$, and define for each $n\geq 1$
 $$
 U_n(g):=\max_{i\in[n]}|\langle g,h_i\rangle|^p\psi_1(\|g\|_2^2),\quad g\in\Li(\xi),
 $$
 and
 $$
 V_n(g):=\varphi(U_n(g))\psi_2(\|g\|_2^2)=\varphi\left(\max_{i\in[n]}|\langle g,h_i\rangle|^p\right)\psi_2(\|g\|_2^2),\quad g\in\Li(\xi).
 $$
Note that $U_n\in\Dom$, $n\geq 1$, by Corollary~\ref{corol_relaxetion_of_cond_on_u} and Lemma~\ref{lemma_estim_for_e_of_max}. Hence, due to Proposition~\ref{prop_chain_rule_for_D}, $V_n$ also belongs to $\Dom$ for all $n\geq 1$.
 
 By the choice of the function $\psi_2$, it is easy to see that for all $g\in L_p^{\uparrow}$
 $$
 V_n(g)\to\varphi(\|g\|_p^p)\psi_2(\|g\|_2^2)=\varphi(\|g\|_p^p),\quad \mbox{as}\ \ n\to\infty,
 $$
  and, consequently, $\{V_n\}_{n\geq 1}$ converges to $\varphi(\|\cdot\|_p^p)$ $\Xi$-a.e., by Corollary~\ref{corol_full_measure_set_for_Xi}. Moreover,
 $$
 |V_n(g)-\varphi(\|g\|_p^p)|\leq 2\|\varphi\|_{\infty}\I_{\{\|g\|_2^2\leq M+1\}},\quad n\geq 1.
 $$
 The dominated convergence theorem implies that $\{V_n\}_{n\geq 1}$ converges to $\varphi(\|\cdot\|_p^p)$ in $\Ll$.
 
 Next, using Proposition~\ref{prop_chain_rule_for_D} and Lemma~\ref{lemma_estim_for_e_of_max}, we can estimate
 \begin{align}
 \e(V_n,V_n)&\leq\frac{1}{2}\|\varphi'\|_{\infty}^2\|\psi_2\|_{\infty}^2\int_{\Li(\xi)}\|\D U_n\|_2^2\Xi(dg)\\
 &+2\|\varphi\|_{\infty}^2\int_{\Li(\xi)}\left(\psi'_2(\|g\|_2^2)\right)^2\|g\|_2^2\Xi(dg)\\
 &\leq\frac{1}{2}\|\varphi'\|_{\infty}^2\|\psi_2\|_{\infty}^2\max_{i\in[n]}\int_{\Li(\xi)}\Big[\psi^2_1(\|g\|_2^2)p^2|\langle g,h_i\rangle|^{2p-2}\|\pr_gh_i\|_2^2\\
 &+4|\langle g,h_i\rangle|^p\left(\psi'_1(\|g\|_2^2)\right)^2\|g\|_2^2\Big]\Xi(dg)\\
 &+2\|\varphi\|_{\infty}^2\|\psi_2'\|_{\infty}^2\int_{\Li(\xi)}\|g\|_2^2\I_{\{\|g\|_2^2\leq M+1\}}\Xi(dg)\\
 &\leq \frac{1}{2}p^2\|\varphi'\|_{\infty}^2\|\psi_2\|_{\infty}^2\|\psi_1\|_{\infty}^2\\
 &\cdot\max_{i\in[n]}\int_{\Li(\xi)}|\langle g,h_i\rangle|^{2p-2}\|\pr_gh_i\|_2^2\I_{\{\|g\|_2^2\leq M+1\}}\Xi(dg)\\
 &+2\|\varphi'\|_{\infty}^2\|\psi_2\|_{\infty}^2\|\psi_1\|_{\infty}^2\int_{\Li(\xi)}|\langle g,h_i\rangle|^p\|g\|_2^2\I_{\{\|g\|_2^2\leq 2M+2\}}\Xi(dg)\\
 &+2\|\varphi\|_{\infty}^2\|\psi_2'\|_{\infty}^2\int_{\Li(\xi)}\|g\|_2^2\I_{\{\|g\|_2^2\leq M+1\}}\Xi(dg).
 \end{align}
Using H\"{o}lder's inequality $|\langle g,h_i\rangle|\leq\|h_i\|_q\|g\|_p=\|g\|_p$ and Lemma~\ref{lemma_estimation_of_int_of_norm_with_peoj}, we have that 
 $$
 \sup_{n\in\N}\e(V_n,V_n)<\infty,
 $$
 if $p\in\left[2,\frac{5}{2}\right]$.
 
 Hence,~\cite[Lemma~I.2.12]{Ma:1992} yields $\varphi(\|\cdot\|_p^p)\in\Dom$ and 
 \begin{equation}\label{f_estimation_of_e_varphi_p}
 \e(\varphi(\|\cdot\|_p^p),\varphi(\|\cdot\|_p^p))\leq\liminf_{n\to\infty}\e(V_n,V_n).
 \end{equation}
 
 In order to compute $\D\varphi(\|\cdot-g_0\|_2^2)$, we take an orthonormal basis $\{h_n\}_{n\geq 1}$ in $L_2$ and note that
 $$
 \|g-g_0\|^2=\sum_{n=1}^{\infty}(\langle g,h_n\rangle-\langle g_0,h_n\rangle)^2.
 $$
 Taking $\psi\in\Cfo(\R)$ such that $\psi=1$ on an interval $[-M,M]$ that contains $\supp\varphi$ and setting 
 $$
 W_n(g)=\varphi\left(\sum_{i=1}^n(\langle g,h_i\rangle-\langle g_0,h_i\rangle)^2\right)\psi(\|g\|_2^2),\quad g\in\Li(\xi),
 $$
a simple calculation shows that 
 $$
 W_n\to\varphi(\|\cdot-g_0\|_2^2)
 $$
 and 
 $$
 \|\D W_n-\D\varphi(\|\cdot-g_0\|_2^2)\|_2\to 0
 $$
 in $\Ll$ as $n\to\infty$. The lemma is proved. 
\end{proof}

\begin{corollary}\label{cor_norm_belongs_to_Dom}
 For each $\varphi\in\Cfo(\R)$ and $g_0\in\Li(\xi)$ the function $U=\|\cdot-g_0\|_2\varphi(\|\cdot\|_2^2)$ belongs to $\Dom$. Moreover, $\|D U\|\leq 1$ $\Xi$-a.e. on $B_r=\{g\in\Li(\xi):\ \|g\|_2\leq r\}$, if $\varphi=1$ on $[-r^2,r^2]$. 
\end{corollary}

\begin{proof}
We take $\psi\in\Cfo(\R)$ such that $\psi=1$ on an interval~$[-M,M]$ that contains $\supp\varphi$.
 For each $\delta>0$, we set
 $$
 V_{\delta}(g)=\left(\|g-g_0\|_2^2\vee\delta^2\right)\psi(\|g\|_2^2),\quad g\in\Li(\xi).
 $$
 Let $\psi_{\delta}\in\Cfb(\R)$ and $\psi_{\delta}(x)=\sqrt{|x|}$ for all $\delta\leq |x|\leq \sup_{g}|V_{\delta}(g)|$. Then by Lemmas~\ref{lemma_estim_for_e_of_max},~\ref{lemma_function_with_compact_supp} and Proposition~\ref{prop_chain_rule_for_D}, the function $U_{\delta}=\psi_{\delta}(V_{\delta})\varphi(\|\cdot\|_2^2)$ belongs to $\Dom$ and 
 $$
 \e(U_{\delta},U_{\delta})\leq C<\infty
 $$
 for all $\delta>0$. Since $U_{\delta}\to U=\|\cdot-g_0\|_2\varphi(\|\cdot\|_2^2)$ in $\Ll$ as $\delta\to 0$, the function $U$ belongs to $\Dom$, by~\cite[Lemma~I.2.12]{Ma:1992}. 
 
 A simple calculation shows that $\|\D U_{\delta}\|\leq 1$ $\Xi$-a.e. on $B_r$  due to the equality $\varphi=1$ on $[-r^2,r^2]$. Hence, by~\cite[Lemma~I.2.12]{Ma:1992}, $\|\D U\|\leq 1$ $\Xi$-a.e. on $B_r$. 
\end{proof}

Let $\FCC$ be the linear span of the set of functions on $\Li(\xi)$ which have a form
\begin{equation}\label{f_form_of_U_with_kompact_supp}
U=u(\langle\cdot,h_1\rangle,\ldots,\langle\cdot,h_m\rangle)\varphi(\|\cdot\|_p^p)=u(\langle\cdot,\vec{h}\rangle)\varphi(\|\cdot\|_p^p),
\end{equation}
where $p\in\left(2,\frac{5}{2}\right]$, $u\in\Cfb(\R^m),$ $\varphi\in\Cfo(\R)$ and $h_j\in L_2(\xi),$  $j\in[m]$.

\begin{remark}
 Each function from $\FCC$ has a compact support in $\Li(\xi)$, by~\cite[Lemma~5.1]{Konarovskyi:2017:EJP}.
\end{remark}

\begin{proposition}\label{prop_density_of_FCC}
 The set $\FCC$ is dense in $\Dom$ with respect to the norm $\e_1^{\frac{1}{2}}$.
\end{proposition}

\begin{proof}
We first note that by Proposition~\ref{prop_chain_rule_for_D} and Lemma~\ref{lemma_function_with_compact_supp}, $\FCC\subset\Dom$. 

To prove the proposition, it is enough to show that each element of $\FC$ can be approximated by elements from $\FCC$. Therefore, let $U\in\FC$ be given by~\eqref{f_form_of_U}, i.e. $U=u(\langle\cdot,\vec{h}\rangle)\varphi(\|\cdot\|_2^2)$. By the dominated convergence theorem and Lemma~\ref{lem_p_norm conv_to_norm},
$$
U_p=u(\langle\cdot,\vec{h}\rangle)\varphi(\|\cdot\|_p^p)\to U\quad\mbox{in}\ \ \Ll\ \ \mbox{as}\ \ p\downarrow 2.
$$
Next, using Proposition~\ref{prop_chain_rule_for_D}, we can estimate, 
\begin{align}
 \e(U_p,U_p)&=\frac{1}{2}\int_{\Li(\xi)}\|\D U_p(g)\|^2_2\Xi(dg)\\
 &\leq 2^{m-1}\sum_{j=1}^m\int_{\Li(\xi)}\varphi^2(\|g\|_p^p)(\partial_ju(\langle g,\vec{h}\rangle))^2\|\pr_g h_j\|^2_2\Xi(dg)\\
 &+2^{m-1}\int_{\Li(\xi)}(u(\langle g,\vec{h}\rangle))^2\|\D \varphi(\|\cdot\|_p^p)(g)\|^2_2\Xi(dg)\\
 &\leq 2^{m-1}\|\varphi\|_{\infty}^2\sum_{j=1}^m\|\partial_ju\|_{\infty}^2\|h_j\|_2^2\int_{\Li(\xi)}\varphi^2(\|g\|_p^p)\Xi(dg)\\
 &+\|u\|_{\infty}^2\e(\varphi(\|\cdot\|_p^p),\varphi(\|\cdot\|_p^p))<C
\end{align}
uniformly in $p\in\left(2,\frac{5}{2}\right]$, by the estimate~\eqref{f_estimation_of_e_varphi_p}, Lemma~\ref{lemma_estimation_of_int_of_norm_with_peoj} and the inequality $\|g\|_2\leq\|g\|_p$ for $p>2$.

Hence, by~\cite[Lemma~I.2.12]{Ma:1992}, there exists a subsequence $\{U_{p_k}\}_{k\geq 1}$ for $p_k\downarrow 2$ such that its Cesaro mean
$$
V_n=\frac{1}{n}\sum_{k=1}^nU_{n_k}\to U
$$
in $\Dom$ (w.r.t. $\e_1^{\frac{1}{2}}$-norm) as $n\to\infty$. Since, $\FCC$ is linear, $V_n\in\FCC$, $n\in\N$. This gives the needed approximation that completes the proof of the proposition.
\end{proof}

\subsection{Quasi-regularity and local property of \texorpdfstring{$(\e,\Dom)$}{(E,D)}}

The aim of this section is to show that $(\e,\Dom)$ is a quasi-regular Dirichlet form.
Let
$$
\Dom_K=\left\{U\in\Dom:\ U=0\ \ \Xi\mbox{-a.e.}\ \mbox{on}\ \Li(\xi)\setminus K\right\}.
$$
We recall that  an increasing sequence $\{K_n\}_{n\geq 1}$ of closed subsets of $\Li(\xi)$ is called an $\e$-\textit{nest}\footnote{The definitions of $\e$-nest, $\e$-quasi-continuity, quasi-regularity and local property are taken from~\cite{Ma:1992} (see Definitions~III.2.1,~III.3.2,~IV.3.1 and~V.1.1, respectively)} if $\bigcup_{n=1}^{\infty}\Dom_{K_n}$ is dense in $\Dom$ (w.r.t. $\e^{\frac{1}{2}}$-norm).

\begin{proposition}
 The Dirichlet form $(\e,\Dom)$ is quasi-regular, that is, it has the following properties
 \begin{enumerate}
  \item[(i)] there exists an $\e$-nest $\{K_n\}_{n\geq 1}$ consisting of compact sets;
  
  \item[(ii)] there exists a dense subset of $\Dom$ (w.r.t. $\e_1^{\frac{1}{2}}$-norm) whose elements have $\e$-quasi-continuous $\Xi$-version;
  
  \item[(iii)] there exist $U_n\in\Dom$, $n\in\N$, having $\e$-quasi-continuous $\Xi$-version $\widetilde{U}_n$, $n\in\N$, and there exists an $\e$-exceptional set $A\subset\Li(\xi)$ such that $\{\widetilde{U}_n,\ n\in\N\}$ separates points of~$\Li(\xi)\setminus A$.
 \end{enumerate}
\end{proposition}

\begin{proof}
 Properties~$(ii)$ and~$(iii)$ follow from the fact that $\FC$ is dense in $\Dom$ (w.r.t. $\e^{\frac{1}{2}}$-norm) and $\FC$ separates points, since $\{\langle\cdot,h\rangle,\ h\in L_2\}$ separates the points of~$\Li(\xi)$.
 
 To prove $(i)$, we set 
 $$
 K_n=\left\{g\in\Li(\xi):\ \|g\|_{2+\frac{1}{n}}\leq n\right\}.
 $$
 Then $\{K_n\}_{n\geq 1}$ is an increasing sequence of compact sets, by~\cite[Lemma~5.1]{Konarovskyi:2017:EJP}. Moreover, it is easily seen that
 $$
 \FCC\subseteq\bigcup_{n=1}^{\infty}\Dom_{K_n}.
 $$
 Consequently, Proposition~\ref{prop_density_of_FCC} yields~$(i)$. It completes the proof of the proposition. 
\end{proof}

\begin{proposition}\label{prop_local_prop}
 The Dirichlet form $(\e,\Dom)$ has the \ local \  property,\  that is,\  ${\e(U,V)=0}$ for all $U,V\in\Dom$ with $\supp(U\cdot\Xi)\cap\supp(V\cdot\Xi)=\emptyset$ and $\supp(U\cdot\Xi)$, $\supp(V\cdot\Xi)$ compact.
\end{proposition}

\begin{proof}
 Let $U\in\Dom$ with $K_U:=\supp(U\cdot\Xi)$ being compact. We first note that the equality $U=0$ $\Xi$-a.e. on a ball $B_r(g_0)=\{g\in\Li(\xi):\ \|g-g_0\|_2<r\}$ implies $\D U=0$ $\Xi$-a.e. on $B_r(g_0)$. Indeed, let $K_U\subset B_R(g_0)$ for some constant $R>0$. We take $\eps\in(0,1)$ and $\varphi\in\Cfo(\R)$ such that $\varphi(x)=0$ for all $|x|\leq (1-\eps)r^2$ and $\varphi(x)=1$ for all $r^2\leq |x|\leq R^2$. Then by Lemmas~\ref{lemma_function_with_compact_supp} and~\ref{lemma_product_belongs_to_Dom}, we can conclude that for all $g\in\Li(\xi)$
 \begin{align}
 \D U(g)&=\D\left[U\varphi(\|\cdot-g_0\|_2^2)\right](g)\\
 &=(\D U(g))\varphi(\|g-g_0\|_2^2)+2U(g)\varphi'(\|g-g_0\|_2^2)g.
 \end{align}
 Hence $\D U(g)=0$ $\Xi$-a.e. on $B_{(1-\eps)r}(g_0)$. Since $\eps$ is arbitrary, we obtain $\D U=0$ $\Xi$-a.e. on $B_r(g_0)$.
 Therefore, the statement easily follows from~\eqref{f_dirichlet_form} and the observation above. The proposition is proved.  
\end{proof}

We also remark that the Dirichlet form $(\e,\Dom)$ satisfies a type of the local property according to the definition from~\cite{Ariyoshi:2005,Bouleau:1991}, that will be needed in Section~\ref{intrinsic metric}.

\begin{lemma}\label{lemma_local_prop_2}
 For each $U\in\Dom$ and $F,G\in C_b^1(\R)$ with $\supp F\cap\supp G=\emptyset$,
 $$
 \e(F(U)-F(0),G(U)-G(0))=0.
 $$
\end{lemma}

\begin{proof}
 The lemma directly follows from Proposition~\ref{prop_chain_rule_for_D}. 
\end{proof}

\subsection{Strictly quasi-regularity and conservativeness in a partial case}

In this section, we will suppose that $\xi$ is constant on some neighbourhoods of~$0$ and~$1$, i.e. there exists $\delta\in\left(0,\frac{1}{2}\right)$ such that $\xi(u)=\xi(0)$, $u\in[0,\delta)$, and $\xi(u)=\xi(1)$, $u\in(1-\delta,1]$. We also set
\begin{equation}\label{f_functions_h_i}
h_1=\frac{1}{\delta}\I_{[0,\delta)}\quad\mbox{and}\quad h_2=\frac{1}{\delta}\I_{[1-\delta,1]}.
\end{equation}
In this case, the space $\Li(\xi)$ is locally compact, that follows from~\cite[Lemma~5.1]{Konarovskyi:2017:EJP} and the following statement.

\begin{lemma}\label{lemma_estim_p_norm_via_2_norm}
For all $p\geq 2$ and $g\in\Li(\xi)$
$\|g\|_p\leq|\langle g,h_1\rangle|\vee|\langle g,h_2\rangle|\leq\frac{1}{\sqrt{\delta}}\|g\|_2$.
\end{lemma}

\begin{proof}
Since $g\in\Li(\xi)$, Proposition~\ref{prop_xi_measurability} implies that $g$ is constant on $[0,\delta)$ and $(1-\delta,1]$. Thus, 
$$
\langle g,h_1\rangle=g(0)\quad\mbox{and}\quad\langle g,h_2\rangle=g(1).
$$
Moreover, $|g(u)|\leq|g(0)|\vee|g(1)|$ for all $u\in(0,1)$, since $g\in D^{\uparrow}$.
Hence, using the Cauchy-Schwarz inequality, we obtain
\begin{align}
\|g\|_p\leq|g(0)|\vee|g(1)|=|\langle g,h_1\rangle|\vee|\langle g,h_2\rangle|\leq\frac{1}{\sqrt{\delta}}\|g\|_2.
\end{align}
The lemma is proved. 
\end{proof}

\begin{proposition}\label{proposition_conservativeteness}
The Dirichlet form $(\e,\Dom)$ is strictly quasi-regular and conservative.
\end{proposition}

\begin{proof}
To prove the strictly quasi-regularity, it is enough to check that $(\e,\Dom)$ is regular\footnote{see e.g. the definition on p.118~\cite{Ma:1992}} according to~\cite[Proposition~V.2.12]{Ma:1992}. Hence, it is needed to prove that $\FC$ is dense in $C_0(\Li(\xi))$ with respect to the uniform norm, where $C_0(\Li(\xi))$ denotes the space of continuous functions on $\Li(\xi)$ with compact support. But this easily follows from the Stone-Weierstrass theorem, Remark~\ref{rem_prop_of_FC} and the fact that each ball in $\Li(\xi)$ is a compact set.

The conservativeness of $(\e,\Dom)$ will follow from~\cite[Theorem~1.6.6]{Fukushima:2011}. Thus, it is enough to show that there exists a sequence $\{U_n,\ n\geq 1\}\subset\Dom$ such that 
\begin{equation}\label{f_consetv_1}
0\leq U_n\leq 1,\quad \lim_{n\to\infty}U_n=1\quad\Xi\mbox{-a.e.}
\end{equation}
and
\begin{equation}\label{f_consetv_2}
\lim_{n\to\infty}\e(U_n,V)=0
\end{equation}
for all $V\in\Dom\cap L_1(\Li(\xi),\Xi)$. 

For each $n\in\N$ we take a function $\psi_n\in\Cfo(\R)$ satisfying 
\begin{itemize}
\item $\supp\psi_n\subset[-2n-1,2n+1]$, $\psi(x)=1$ on $[-n,n]$ and $\psi_n(x)\in[0,1]$ for $n<|x|<2n+1$;

\item $|\psi_n'(x)|\leq\frac{1}{n}$ and $|\psi_n''(x)|\leq\frac{C}{n}$ for all $x\in\R$ and a constant $C$ that is independent of~$n$.
\end{itemize}
We also set
$$
U_n(g)=u_n(\langle g,h_1\rangle,\langle g,h_2\rangle),\quad g\in\Li(\xi)\ \ \mbox{and}\ \ n\geq 1,
$$
where $u_n(x,y)=\psi_n(x)\psi_n(y)$, $x,y\in\R$, and $h_1$, $h_2$ are defined by~\eqref{f_functions_h_i}. Then, by Lemma~\ref{lemma_estim_p_norm_via_2_norm}, for each $\varphi\in\Cfo(\R)$ satisfying $\varphi=1$ on $[-(2n+1)^2,(2n+1)^2]$ the equality
$$
U_n(g)=U_n(g)\varphi(\|g\|_2^2),\quad g\in\Li(\xi),
$$
holds. This implies that $U_n\in\FC$ and
\begin{align}
LU&=\frac{1}{2}\sum_{i,j=1}^2\partial_i\partial_ju_n(\langle g,h_1\rangle,\langle g,h_2\rangle)\langle\pr_gh_i,\pr_gh_j\rangle\\
&+\frac{1}{2}\sum_{j=1}^2\partial_ju_n(\langle g,h_1\rangle,\langle g,h_2\rangle)\langle\xi-\pr_g\xi,h_j\rangle,\quad g\in\Li(\xi),
\end{align}
for all $n\geq 1$, where $L$ is defined by~\eqref{f_generator}. By the construction of $U_n$, $\{U_n,\ n\geq 1\}$ satisfies~\eqref{f_consetv_1}. Moreover, using the Cauchy-Schwarz inequality, the inequality $\|\pr_gh\|_2\leq\|h\|_2$ and the dominated convergence theorem, we have for every $V\in\Dom\cap L_1(\Li(\xi),\Xi)$
\begin{align}
\e(U_n,V)&=-(LU_n,V)_{\Li(\xi)}\\
&=\frac{1}{2}\sum_{i,j=1}^2\int_{\Li(\xi)}\partial_i\partial_ju_n(\langle g,h_1\rangle,\langle g,h_2\rangle)\langle\pr_gh_i,\pr_gh_j\rangle V(g)\Xi(dg)\\
&+\frac{1}{2}\sum_{j=1}^2\int_{\Li(\xi)}\partial_ju_n(\langle g,h_1\rangle,\langle g,h_2\rangle)\langle\xi-\pr_g\xi,h_j\rangle V(g)\Xi(dg)\to 0
\end{align} 
as $n\to\infty$. The proposition is proved. 
\end{proof}

\section{Intrinsic metric associated to \texorpdfstring{$(\e,\Dom)$}{(E,D)}}\label{intrinsic metric}

The aim of this section is to show that $L_2$-metric is the intrinsic metric associated to $(\e,\Dom)$ and to prove the analog of Varadhan's formula. We will use the result obtained in~\cite{Ariyoshi:2005} for the proof of Varadhan's formula (see also~\cite{Masanori:2003} for the Dirichlet forms on $L_2(\mu)$, where $\mu$ is a probability measure).

\subsection{The boundedness of \texorpdfstring{$\D U$}{DU} implies the Lipschitz continuity of \texorpdfstring{$U$}{U}}

In this section we will prove that any function $U$ from $\Dom$ with $\|\D U\|\leq 1$ $\Xi$-a.e. is 1-Lipschitz continuous.  

\begin{proposition}\label{prop_lipschitz_continuity}
 Let $U\in\Dom$ and $\|\D U\|_2\leq 1$ $\Xi$-a.e. on a convex open set $B\subseteq\Li(\xi)$. Then $U$ has an 1-Lipschitz modification $\widetilde{U}$ on $B$, i.e. there exists a function $\widetilde{U}:B\to\R$ such that $\Xi\{g\in B:\ \widetilde{U}(g)\neq U(g)\}=0$ and 
 \begin{equation}\label{f_lipschitz_inequality}
 |\widetilde{U}(g_1)-\widetilde{U}(g_0)|\leq \|g_1-g_0\|_2
 \end{equation}
 for all $g_0,g_1\in B$.
\end{proposition}

\begin{remark}
If $U\in\FC$, then  
$$
U(g_1)-U(g_0)=\int_0^1\langle\D U(g_t),g_1-g_0\rangle dt
$$
for all $g_0,g_1\in\St$, where $g_t=g_0+t(g_1-g_0)$. This follows from the fact that  $\sigma^{\star}(g_t)\supseteq\sigma^{\star}(g_1-g_0)$ for all $t\in(0,1)$ and $g_0,g_1\in\St$. Therefore, the statement holds for all $U\in\FC$. 
\end{remark}

\begin{proof}[Proof of Proposition~\ref{prop_lipschitz_continuity}]
{\it Step I.} First we show that for each $n\geq 1$,~\eqref{f_lipschitz_inequality} holds $\Xi_n$-a.e on $B$. Let $n\geq 2$ be fixed.
Since $\FC$ is dense in $\Dom$ (w.r.t. $\e^{\frac{1}{2}}$-norm), there exists a sequence $\{U_k\}_{k\geq 1}\subset\FC$ such that $U_k\to U$ and $\|\D U_k-\D U\|_2\to 0$ in $\LL$ as $k\to\infty$. Hence, $U_k\to U$ and $\|\D U_k-\D U\|_2\to 0$ in $L_2(\Li(\xi),\Xi_n)$.
 
Let $A\subseteq B$ such that $\Xi(B\setminus A)=0$ and $\|\D U(g)\|\leq 1$ for all $g\in A$. We set
$$
A_n=A\cap\{\chi_n(q,x):\ q\in Q^n,\ x\in E_0^n\}.
$$
Then by Remark~\ref{rem_repres_of_step_functions} and Lemma~\ref{lem_prop_of_Xi_n}~(iii), $\Xi_n(B\setminus A_n)=0$. Since $\Xi_n$ is the push forward of the measure $\mu_{\xi}^n\otimes\lambda_n$ under the map $\chi_n$ (see Lemma~\ref{lem_prop_of_Xi_n}~(i)), it is easy to see that there exists $Q_1\subseteq Q^n$ such that $\mu_{\xi}^n(Q^n\setminus Q_1)=0$ and $\lambda_n(B(q)\setminus A_n(q))=0$ for all $q\in Q_1$, where $A_n(q)=\{x\in E_0^n:\ \chi_n(q,x)\in A_n\}$ and $B(q)=\{x\in E_0^n:\ \chi_n(q,x)\in B\}$.

We next note that
\begin{align}
&\int_{\Li(\xi)}|U_k(g)-U(g)|^2\Xi_n(dg)\\
&=\int_{Q^n}\left[\int_{E^n}|U_k(\chi_n(q,x))-U(\chi_n(q,x))|^2\lambda_n(dx)\right]\mu_{\xi}^n(dq)\to 0
\end{align}
and, similarly,
$$
\int_{Q^n}\left[\int_{E^n}\|\D U_k(\chi_n(q,x))-\D U(\chi_n(q,x))\|_2^2\lambda_n(dx)\right]\mu_{\xi}^n(dq)\to 0
$$
as $k\to\infty$. Consequently, we can choose a subsequence $\{k'\}\subseteq\N$ (we assume that $\{k'\}$ coincides with $\N$ without loss of generality) and a measurable subset $Q_2\subseteq Q^n$ such that  $\mu_{\xi}^n(Q^n\setminus Q_2)=0$ and for all $q\in Q_2$
\begin{align}\label{f_conv_of_deriv_in_E_n}
\begin{split}
\int_{E^n}|U_k(\chi_n(q,x))-U(\chi_n(q,x))|^2\lambda_n(dx)&\to 0,\\
\int_{E^n}\|\D U_k(\chi_n(q,x))-\D U(\chi_n(q,x))\|_2^2\lambda_n(dx)&\to 0
\end{split}
\end{align}
as $k\to \infty$.
 
Let $q\in Q_1\cap Q_2$ be fixed and
\begin{align}
f_k(x):=&U_k(\chi_n(q,x)),\quad x\in E_0^n,\\
f(x):=&U(\chi_n(q,x)),\quad x\in E_0^n.
\end{align}
Then $f_k$, $k\geq 1$, belong to $\Cfo(E^n)$ and 
\begin{equation}\label{f_deriv_of_h_k}
\D U_k(\chi_n(q,x))=\sum_{i=1}^n\frac{\partial f_k(x)}{\partial x_i}\frac{\I_{[q_{i-1},q_i)}}{q_i-q_{i-1}},\quad x\in E_0^n, 
\end{equation}
by Lemma~\ref{lemma_connection_D_with_D_n}. We are going to show that $\D U(\chi_n(q,\cdot))$ is also given by~\eqref{f_deriv_of_h_k}, where the partial derivatives of $f_k$ is replaced by the Sobolev partial derivatives of $f$. 

We note that $\D U(\chi_n(q,\cdot))$ can be given as follows
$$
\D U(\chi_n(q,x))=\sum_{i=1}^n\widetilde{f}^i(x)\frac{\I_{[q_{i-1},q_i)}}{q_i-q_{i-1}},\quad x\in E_0^n,
$$
for some measurable functions $\widetilde{f}^i:E_0^n\to \R$, since the set $\left\{\sum_{i=1}^nx_i\I_{[q_{i-1},q_i)},\  x\in \R^n\right\}$ is closed in $L_2(\xi)$.
Moreover, by~\eqref{f_conv_of_deriv_in_E_n}, we have that
$$
\int_{E_0^n}|f_k(x)-f(x)|^2\lambda_n(dx)\to 0
$$
and
$$
\int_{E_0^n}\sum_{i=1}^n\left[\widetilde{f}^i(x)-\frac{\partial f_k(x)}{\partial x_i}\right]^2(q_i-q_{i-1})\lambda_n(dx)\to 0
$$
as $k\to\infty$. It immediately implies that $f$ belongs to the Sobolev space $H^{1,2}(E_0^n)$ with $\widetilde{f}^i=\frac{\partial f}{\partial x_i}$. In particular,
\begin{equation}\label{f_int_by_parts_f}
\int_{\R^n}f(x)\frac{\partial \varphi(x)}{\partial x_i}dx=-\int_{\R^n}\widetilde{f}^i(x)\varphi(x)dx.
\end{equation}
for each $\varphi\in\Cfo(\R^n)$  with $\supp\varphi\subset E^n_0$ and $f,\ \widetilde{f}^i$, $i\in[n]$, equal zero outside $E^n$.

Next, let $\varphi\in\Cfo(\R^n)$ be a non negative function with
$$
\int_{\R^n}\varphi(x)dx=1.
$$
Then the convolution
$$
f_{\eps}(x)=f\ast\varphi_{\eps}(x)=\int_{\R^n}f(y)\varphi_{\eps}(x-y)dy,\quad x\in\R^n,
$$
where $\varphi_{\eps}(x)=\eps^{-n}\varphi(x\eps^{-1})$, belongs to $C^{\infty}(\R^n)$ and converges to $f$ $\lambda_n$-a.e. on $E_0^n$. Moreover, by~\eqref{f_int_by_parts_f},
$$
\frac{\partial f_{\eps}(x)}{\partial x_i}=\widetilde{f}^i\ast\varphi_{\eps}(x)
$$
for every $x\in E_0^n$ and all $\eps>0$ satisfying $\supp\varphi_{\eps}(x-\cdot)\subset E_0^n$.

We recall that $B(q)=\{x\in E_0^n:\ \chi_n(q,x)\in B\}$. Let $B(q)\neq\emptyset$. It is easily seen that $B(q)$ is an open convex subset of $E_0^n$. Then for each $x\in B(q)$ and $\eps>0$ such that $\supp\varphi_{\eps}(x-\cdot)\subset B(q)$ we can estimate
\begin{align}\label{f_estim_of_deriv_of_f}
\begin{split}
 \sum_{i=1}^n&\left(\frac{\partial f_{\eps}(x)}{\partial x_i}\right)^2\frac{1}{q_i-q_{i-1}}=\sum_{i=1}^n\left(\widetilde{f}^i\ast\varphi_{\eps}(x)\right)^2\frac{1}{q_i-q_{i-1}}\\
 &\leq\sum_{i=1}^n\int_{\R^n}(\widetilde{f}^i(y))^2\varphi_{\eps}(x-y)dy\frac{1}{q_i-q_{i-1}}=\int_{\R^n}\sum_{i=1}^n\frac{(\widetilde{f}^i(y))^2}{q_i-q_{i-1}}\varphi_{\eps}(x-y)dy\\
 &=\int_{E_0^n}\|\D U(\chi_n(q,y))\|_2^2\varphi_{\eps}(x-y)\lambda_n(dy)\\
 &=\int_{B(q)}\|\D U(\chi_n(q,y))\|_2^2\varphi_{\eps}(x-y)\lambda_n(dy)\leq 1,
\end{split}
\end{align}
since $\|\D U(\chi_n(q,\cdot))\|_2\leq 1$ $\lambda_n$-a.e. on $B(q)$.

Let $x^0,x^1\in B(q)$ and $\eps_0>0$ such that $f_{\eps}(x^i)\to f(x^i)$ and $\supp\varphi_{\eps_0}(x^i-\cdot)\subset B(q)$, $i=0,1$. Using the convexity of $B(q)$, it is easy to see that
$$
\supp\varphi_{\eps_0}(x^t-\cdot)\subset B(q), \quad t\in(0,1),
$$
where $x^t=x^0+t(x^1-x^0)$. By H\"{o}lder's inequality and~\eqref{f_estim_of_deriv_of_f}, we can estimate
\begin{align}
 (f_{\eps}(x^1)-f_{\eps}(x^0))^2&=\left(\int_0^1\frac{d}{dt}f_{\eps}(x^t)dt\right)^2=\left(\int_0^1\sum_{i=1}^n\partial_i f_{\eps}(x^t)(x_i^1-x_i^0)dt\right)^2\\
 &\leq\int_0^1\sum_{i=1}^n\left(\partial_i f_{\eps}(x^t)\right)^2\frac{1}{q_i-q_{i-1}}dt\sum_{i=1}^n(x_i^1-x_i^0)^2(q_i-q_{i-1})\\
 &\leq \|\chi_n(q,x^1)-\chi_n(q,x^0)\|_2^2
\end{align}
for all $\eps\in(0,\eps_0]$. Hence using the convergence of $f_{\eps}(x^i)$ to $f(x^i)$, $i=0,1$, and the previous estimate, we have that
\begin{equation}\label{f_lipshitz_ineq}
|U(\chi(q,x^1))-U(\chi(q,x^0))|\leq \|\chi_n(q,x^1)-\chi_n(q,x^0)\|_2.
\end{equation}

Since~\eqref{f_lipshitz_ineq} holds for all $q\in Q_1\cap Q_2$ and $x^i\in B(q)$, $i=0,1$, such that $f_{\eps}(x^i)\to f(x^i)$ as $\eps\to 0$, we have that  
\begin{equation}\label{f_lip_prop_for_st}
|U(g_1)-U(g_0)|\leq \|g_1-g_0\|_2,\quad\Xi_n\mbox{-a.e. on}\ \ B,
\end{equation}
due to the equalities $\mu_{\xi}^n(Q^n\setminus(Q_1\cap Q_2))=0$ and $\lambda_n\{x\in B(q):f_{\eps}(x)\not\to f(x)\}=0$.

We also note that using the same argument, we can show that~\eqref{f_lipshitz_ineq} holds $\Xi_1$-a.e. on~$B$.

\vspace{8pt}
{\it Step II.} Let $\widetilde{B}_n\subseteq B\cap\supp\Xi_n$ such that $\Xi_n(B\setminus\widetilde{B}_n)=0$ and for all $g_0,g_1\in\widetilde{B}_n$ the inequality~\eqref{f_lip_prop_for_st} holds. Since $\Xi_n(B\setminus\widetilde{B}_n)=0$, $\widetilde{B}_n$ is dense in $B\cap\supp\Xi_n$. Consequently, there exists a unique 1-Lipschitz function $\widetilde{U}_n:B\cap\supp\Xi_n\to\R$ that is the extension of $U$ to $B\cap\supp\Xi_n$. Moreover, $\widetilde{U}_n=U$ $\Xi_n$-a.e. By the uniqueness of the extension and  Corollary~\ref{corol_supp_of_Xi_n}, we have that $\widetilde{U}_n=\widetilde{U}_{n+1}$ on $B\cap\supp\Xi_n=B\cap\{g\in\Li:\ \sharp g\leq n\}$. Therefore, we can define 
$$
\widetilde{U}_{\infty}(g)=\widetilde{U}_n(g),\quad g\in B\cap\supp\Xi_n=B\cap\{g\in\Li:\ \sharp g\leq n\}.
$$
Thus, $\widetilde{U}_{\infty}$ is an 1-Lipschitz function defined on $B\cap\left(\bigcup_{n=1}\supp\Xi_n\right)=B\cap\St$, since for any $g_0,g_1\in B\cap\St$ there exists $n\in\N$ such that $g_0,g_1\in B\cap\{g\in\Li:\ \sharp g\leq n\}$. By the density of $B\cap\St$ in $B$, we can extend $\widetilde{U}_{\infty}$ to an 1-Lipschitz function $\widetilde{U}$ defined on $B$. Moreover, $\widetilde{U}=U$ $\Xi$-a.e. on $B$ because $\Xi(\Li(\xi)\setminus\St)=0$, by Corollary~\ref{corol_full_measure_set_for_Xi}, that completes the proof of the proposition.
\end{proof}

\subsection{Intrinsic metric and Varadhan's formula}

Since the measure $\Xi$ is $\sigma$-finite, we will define the intrinsic metric associated to $(\e,\Dom)$ using a localization of the domain $\Dom$ (see~\cite{Ariyoshi:2005}). Let $L_0(\Xi)$ denote the set of all measurable functions on $\Li(\xi)$ and $K_n:=\{g\in\Li(\xi):\ \|g\|_2\leq n\}$, $n\in\N$. Then the family of balls $\{K_n\}_{n\geq 1}$ satisfies the following conditions
\begin{enumerate}
 \item[$(N1)$] For every $n\in\N$ there exists $V_n\in\Dom$ such that $V_n\geq 1$ $\Xi$-a.e. on $K_n$;
 
 \item[$(N2)$] $\bigcup_{n=1}^{\infty}\Dom_{K_n}$ is dense in $\Dom$ (w.r.t. $\e^{\frac{1}{2}}$-norm).
\end{enumerate}

\begin{remark}
 We note that the family $\{K_n\}_{n\geq 1}$ is a nest. It is also a nest according the definition given in~\cite{Ariyoshi:2005}, where the topology (on $\Li(\xi)$) is not needed. 
\end{remark}

We set
$$
\Dom_{loc}(\{K_n\})=\left\{U\in L_0(\Xi):\ \begin{array}{l}
                                            \mbox{there exists}\ \ \{U_n\}_{n\geq 1}\subset\Dom\ \ \mbox{such that}\\
                                            U=U_n\ \ \Xi\mbox{-a.e. on}\ \ K_n\ \ \mbox{for each}\ \ n
                                           \end{array}
\right\}
$$
and let $\Dom_{loc,b}(\{K_n\})$ denote the set of all essentially bounded functions from $\Dom_{loc}(\{K_n\})$. For $U,V\in\Dom_b$, where $\Dom_b$ is the set of all essentially bounded functions from $\Dom$, we define
$$
I_U(V)=2\e(UV,U)-\e(U^2,V).
$$
By the locality of $(\e,\Dom)$ (see Lemma~\ref{lemma_local_prop_2}), $I_U(V)$ and $\D U$ can be well-defined for all $U\in\Dom_{loc,b}(\{K_n\})$ and $V\in\bigcup_{n=1}^{\infty}\Dom_{K_n,b}$, where $\Dom_{K_n,b}=\Dom_{K_n}\cap\Dom_b$, setting $I_U(V)=I_{U_n}(V)$ and $\D U=\D U_n$ if $V\in\Dom_{K_n,b}$ and $U_n=U$ $\Xi$-a.e. on $K_n$. 

We set 
$$
\Dom_0=\left\{U\in\Dom_{loc,b}(\{K_n\}):\ \ I_U(V)\leq\|V\|_{L_1(\Xi)}\ \ \mbox{for every}\ \ V\in\bigcup_{n=1}^{\infty}\Dom_{K_n,b}\right\}.
$$

\begin{remark}
 According to~\cite[Proposition~3.9]{Ariyoshi:2005}, the set $\Dom_0$ does not depend on the family of increasing sets $\{K_n\}_{n\geq 1}$ that satisfies $(N1)$, $(N2)$.
\end{remark}

\begin{lemma}\label{lemma_set_Dom_0}
 The set $\Dom_0$ coincides with $\left\{U\in\Dom_{loc,b}(\{K_n\}):\ \ \|\D U\|_2\leq 1\ \ \Xi\mbox{-a.e.}\right\}$. 
\end{lemma}

\begin{proof}
 The statement easily follows from the relation 
 $$
 I_U(V)=\int_{\Li(\xi)}\|\D U(g)\|_2^2V(g)\Xi(dg),
 $$ 
 the density of $\FC_{K_n}=\{U\in\FC:\ U=0\ \ \Xi\mbox{-a.e. on}\ \ \Li(\xi)\setminus K_n\}$ in $L_1(K_n,\Xi)$ (w.r.t. $L_1$-norm) and the duality between $L_1(K_n,\Xi)$ and $L_{\infty}(K_n,\Xi)$. 
\end{proof}

We note that each $U\in\Dom_0$ has a continuous modification, by Lemma~\ref{lemma_set_Dom_0} and Proposition~\ref{prop_lipschitz_continuity}. Therefore, considering such a function, we will take its continuous modification.

\begin{theorem}\label{theorem_intrinsic_metric}
 The intrinsic metric for the Dirichlet form $(\e,\Dom)$ is the $L_2$-metric, that is, for all $g_0,g_1\in\Li(\xi)$
 \begin{align}
 \|g_1-g_0\|_2&=\sup_{U\in\Dom_0}\{U(g_1)-U(g_0)\}\\
 &=\sup\left\{U(g_1)-U(g_0):\ \ U\in\Dom_{loc,b}(\{K_n\}),\ \ \|\D U\|_2\leq 1\ \ \Xi\mbox{-a.e.}\right\}.
 \end{align}
\end{theorem}

\begin{proof}
 The equality
 \begin{align}
 \sup_{U\in\Dom_0}\{U(g_1)-U(g_0)\}=\sup\left\{U(g_1)-U(g_0):\ \ 
 \begin{array}{l} 
 U\in\Dom_{loc,b}(\{K_n\}),\ \ 
 |\D U\|_2\leq 1\ \ \Xi\mbox{-a.e.}
 \end{array}\right\}
 \end{align}
 follows from Lemma~\ref{lemma_set_Dom_0}. Proposition~\ref{prop_lipschitz_continuity} implies the lower bound 
 \begin{align}
 \|g_1-g_0\|\geq\sup\left\{U(g_1)-U(g_0):\ \ 
 \begin{array}{l} 
 U\in \Dom_{loc,b}(\{K_n\}),\ \ 
 \|\D U\|_2\leq 1\ \ \Xi\mbox{-a.e.}
 \end{array}\right\}.
 \end{align}
 To finish the proof, for $g_0,g_1\in\Li(\xi)$ and $g_0\neq g_1$ we need to find $U\in\Dom_0$ such that $U(g_1)-U(g_0)=\|g_1-g_0\|_2$. We take $u\in C_b^1(\R)$ such that $u(x)=x$ for all $|x|\leq \|g_1\|_2\vee\|g_0\|_2$ and $|u'(x)|\leq 1$, $x\in\R$, and define
 $$
 U(g)=u\left(\frac{\langle g,g_1-g_0\rangle}{\|g_1-g_0\|_2}\right),\quad g\in\Li(\xi).
 $$
 Since $\frac{|\langle g_i,g_1-g_0\rangle|}{\|g_1-g_0\|_2}\leq\|g_0\|_2\vee\|g_1\|_2$, we have
 $$
 U(g_1)-U(g_0)=\|g_1-g_0\|_2.
 $$
 Moreover, it is easy to see that $U\in\Dom_{loc,b}$ and
 $$
 \D U(g)=u'\left(\frac{\langle g,g_1-g_0\rangle}{\|g_1-g_0\|_2}\right)\frac{\pr_g(g_1-g_0)}{\|g_1-g_0\|_2},
 $$
 by Proposition~\ref{prop_chain_rule_for_D}. Consequently, $\|\D U(g)\|_2\leq 1$ for all $g\in\Li(\xi)$. This completes the proof of the theorem. 
\end{proof}

Next, let $\{T_t\}_{t\geq 0}$ denote the semigroup on $\LL$ associated with $(\e,\Dom)$. For measurable sets $A,B\subseteq\Li(\xi)$ with positive $\Xi$-measure we define
$$
P_t(A,B)=\int_{\Li(\xi)}\I_A(g)\cdot T_t\I_B(g)\Xi(dg)
$$
and 
$$
\dm(A,B)=\essinf\{\|g-f\|_2:\ g\in A,\ f\in B\}.
$$

\begin{theorem}\label{theorem_varadhan_formula}
 For any measurable $A,B\subset\Li(\xi)$ with $0<\Xi(A)<\infty$, $0<\Xi(B)<\infty$ and $A$ or $B$ open the relation
 $$
 \lim_{t\to 0}t\ln P_t(A,B)=-\frac{\dm(A,B)^2}{2}
 $$
 holds.
\end{theorem}

\begin{proof}
 The statement follows from the general result for symmetric
diffusions obtained in~\cite[Theorem~2.7]{Ariyoshi:2005} and Theorem~\ref{theorem_intrinsic_metric}. 
\end{proof}

The following result is a consequence of Theorem~5.2~\cite{Ariyoshi:2005} and Theorem~\ref{theorem_intrinsic_metric}. Let $\|g-A\|_2:=\essinf_{f\in A}\|g-f\|_2$, \ $g\in\Li(\xi)$.

\begin{theorem}\label{theorem_varadhan_formula_spec_case}
 Let $A$ be a non empty open 
 subset of $\Li(\xi)$ with $\Xi(A)<\infty$ and $\Theta$ be any probability measure which is mutually absolutely continuous with respect to $\Xi$. Then the function $u_t=-t\ln T_t\I_A$ converges to $\frac{\|\cdot-A\|_2^2}{2}$ in the following senses.
 \begin{enumerate}
  \item[(a)] $u_t\cdot\I_{\{u_t<\infty\}}$ converges to $\frac{\|\cdot-A\|_2^2}{2}\cdot\I_{\{\|\cdot-A\|_2<\infty\}}$ in $\Theta$-probability as $t\to 0$.
  
  \item[(b)] If $F$ is a bounded function on $[0,\infty]$ that is continuous on $[0,\infty)$, then $F(u_t)$ converges to $F\left(\frac{\|\cdot-A\|_2^2}{2}\right)$ in $L_2(\Li(\xi),\Theta)$ as $t\to 0$.
 \end{enumerate}
\end{theorem}

\section{Sticky-reflected particle system}\label{section_process}

In this section, we will study some properties of the process associated with the Dirichlet form $(\e,\Dom)$. Let $X=\left(\Omega,\F,(\F_t)_{t\geq 0},\{X_t\}_{t\geq 0}, \{\p_g\}_{g\in\Li(\xi)_{\Delta}}\right)$ be a $\Xi$-tight (Markov) diffusion\footnote{see~\cite[Definition~V.1.10]{Ma:1992}} process with state space $\Li(\xi)$ and life time $\zeta$ that is properly associated with $(\e,\Dom)$. Such a process $X$ exists and is unique up to $\Xi$-equivalence according to~\cite[Theorems~IV.6.4 and~V.1.11]{Ma:1992}. We recall that 
$X$ is continuous on $[0,\zeta)$, i.e.
$$
\p_g\left\{t\mapsto X_t\ \ \mbox{is continuous on}\ \ [0,\zeta)\right\}=1\quad \mbox{for}\ \ \e\mbox{-q.e.}\ \ g\in\Li(\xi).
$$

We also remark that by Proposition~\ref{proposition_conservativeteness}, $\p_g\{\zeta<\infty\}=0$ for $\e$-q.e. $g\in\Li(\xi)$, if $\xi$ is constant on some neighbourhoods of $0$ and $1$.

\subsection{\texorpdfstring{$X$}{X} as \texorpdfstring{$L_2(\xi)$}{L2(xi)}-valued semimartingale}

In this section, we will show that the process $X_t$, $t\in[0,\zeta)$, is a continuous local semimartingale in $\Li(\xi)$ under $\p_g$ for $\e$-q.e. $g\in\Li(\xi)$. Letting
$$
M_t=X_t-\frac{1}{2}\int_0^t(\xi-\pr_{X_s}\xi)ds,\quad t\in[0,\zeta),
$$
the following theorem holds.

\begin{theorem}\label{theorem_martingale_problem_for_X}
There exists an $\e$-exceptional subset $N$ of $\Li(\xi)$ such that for all $g\in\Li(\xi)\setminus N$ and each $(\F_t)$-stopping time $\tau$ satisfying $\p_{g}\{\tau<\zeta\}=1$ and $\E_g\|X_t^{\tau}\|_2^2<\infty$, $t\geq 0$, the process $M_t^{\tau}$, $t\geq 0$, is a continuous square integrable $(\F_t)$-martingale under $\p_g$ in $L_2(\xi)$ with the  quadratic variation\footnote{see~\cite[Definition 2.9]{Gawarecki:2011} for the precise definition of quadratic variation of Hilbert-space-valued martingales}
$$
[M^{\tau}_{\cdot}]_t=\int_0^{t\wedge\tau}\pr_{X_s}ds,\quad t\geq 0,
$$ 
where $X_t^{\tau}:=X_{t\wedge\tau}$ and $M_t^{\tau}:=M_{t\wedge\tau}$. In particular, for each $h_1,h_2\in L_2(\xi)$ the processes $\langle M_t^{\tau},h_i\rangle$, $t\geq 0$, $i\in[2]$, are continuous square integrable $(\F_t)$-martingales under $\p_g$ with the joint quadratic variation
$$
[\langle M^{\tau}_\cdot,h_1\rangle,\langle M^{\tau}_\cdot,h_2\rangle]_t=\int_0^{t\wedge\tau}\langle\pr_{X_s}h_1,h_2\rangle ds,\quad t\geq 0.
$$
\end{theorem}

\begin{proof}
 The statement easily follows from the martingale problem for $X$ (see, e.g.,~\cite[Theorem~3.4~$(i)$]{MR1335494}) and the fact that for all $\varphi\in\Cfo(\R)$ with $\varphi=1$ on an interval $[-C,C]$ and $U(g):=\langle g,h\rangle\varphi(\|g\|_2^2)$, $g\in\Li(\xi)$, we have 
$$
\D U(g)=\pr_gh\quad\mbox{and}\quad LU(g)=\frac{1}{2}\langle\xi-\pr_{g}\xi,h\rangle
$$ 
for all $g\in\Li(\xi)$  satisfying $\|g\|_2^2\leq C$. 
\end{proof}

\begin{corollary}\label{corollary_M_is_martingale}
If $\xi$ is a constant on some neighbourhoods of $0$ and $1$, then for $\e$-q.e. $g\in\Li(\xi)$ $\E_g\|X_t\|_2^2<\infty$, $t\geq 0$, and the process $M_t$, $t\geq 0$, is a continuous square integrable $(\F_t)$-martingale under $\p_g$ in $L_2(\xi)$ with the quadratic variation
$$
[M_{\cdot}]_t=\int_0^t\pr_{X_s}ds,\quad t\geq 0.
$$ 
\end{corollary}

\begin{proof}
The statement of the corollary follows from Theorem~\ref{theorem_martingale_problem_for_X} and Proposition~\ref{proposition_conservativeteness}.
\end{proof}

\subsection{Evolution of the empirical mass process}

Let $\Pp_2$ denote the space of probability measures on $\R$ with the finite second moment. We recall that $\Pp_2$ is a Polish space with respect to the quadratic Wasserstein metric 
\begin{equation}\label{equ_wasserstein_distance}
d_{\mathcal W}(\nu_1,\nu_2)=\left(\inf_{\nu\in\chi(\nu_1,\nu_2)}\iint_{\R^2}|x-y|^2\nu(dx,dy)\right)^{\frac{1}{2}},
\end{equation}
where $\chi(\nu_1,\nu_2)$ denotes the set of all probability measures on $\R^2$ with marginals $\nu_1,\nu_2\in\Pp_2$. Let $\iota g$ denote the push forward of the Lebesgue measure $\leb$ on $[0,1]$ under $g\in\Li(\xi)$, that is,
$$
\iota g(A)=\leb\{u:\ g(u)\in A\}, \quad A\in\B(\R).
$$

\begin{remark}
The map $\iota$ is  bijective isometry between $L_2$ and $\Pp_2$ (for more details see, e.g.,~\cite[Section~2.1]{Brenier:2013}).
\end{remark} 

Let 
\begin{equation}\label{f_mu}
\mu_t:=\iota X(\cdot,t),\quad t\geq 0,
\end{equation}
where $\iota\Delta:=\Delta$.
We are going to show that the process $\mu_t$, $t\geq 0$, is a martingale solution on $[0,\zeta)$ to the stochastic partial differential equation 
\begin{equation}\label{f_SPDE_for_mu}
d\mu_t =  \Gamma(\mu_t)dt + \mathop{{\rm div}}( \sqrt {\mu_t} dW_t),
\end{equation}
with $\langle\alpha,\Gamma(\nu)\rangle=\frac{1}{2}\sum_{x\in\supp\nu}\alpha''(x)$,\ \ $\alpha\in\Cfo(\R)$. In particular, it will yield that~\eqref{f_SPDE_for_mu} has no unique solution, since the modified massive Arratia flow is a martingale solution to the same equation (see~\cite[Section~1.3.1]{Konarovskyi_LDP:2015}).

\begin{proposition}\label{prop_int_alpha_to_dom}
 For each $\alpha\in C^1_b(\R)$ and $\varphi\in\Cfo(\R)$ the function 
 $$
 U(g)=\int_0^1\alpha(g(s))ds\cdot\varphi(\|g\|_2^2),\quad g\in\Li(\xi),
 $$
 belongs to $\Dom$ and 
 $$
 \D U(g)=\alpha'(g)\varphi(\|g\|_2^2)+\int_0^1\alpha(g(s))ds\cdot 2\varphi'(\|g\|_2^2)g,\quad g\in\Li(\xi).
 $$
\end{proposition}

\begin{proof}
 The proof is given in the appendix. 
\end{proof}

\begin{corollary}\label{corol_int_alpha_to_dom}
 Let $\alpha_j\in C^1_b(\R)$, $j\in[m]$, $\varphi\in\Cfo(\R)$ and $u\in C_b^1(\R^m)$. Then the function
 \begin{equation}\label{f_U_of_int}
 \begin{split}
 U(g)&=u\left(\int_0^1\alpha_1(g(s))ds,\ldots,\int_0^1\alpha_m(g(s))ds\right)\varphi(\|g\|_2^2)\\
 &=u\left(\int_0^1\vec{\alpha}(g(s))ds\right)\varphi(\|g\|_2^2),\quad g\in\Li(\xi),
 \end{split}
 \end{equation}
 belongs to $\Dom$ and 
 \begin{align}\label{f_deriv_of_u_of_int}
 \begin{split}
 \D U(g)&=\sum_{j=1}^m\partial_ju\left(\int_0^1\vec{\alpha}(g(s))ds\right)\alpha'(g)\varphi(\|g\|_2^2)\\
 &+u\left(\int_0^1\vec{\alpha}(g(s))ds\right)\cdot 2\varphi'(\|g\|_2^2)g,\quad g\in\Li(\xi).
 \end{split}
 \end{align}
\end{corollary}

\begin{proof}
 The corollary follows from Propositions~\ref{prop_chain_rule_for_D} and~\ref{prop_int_alpha_to_dom}.
\end{proof}

\begin{proposition}\label{prop_integr_by_parts_for_mu}
 Let $\alpha_j\in C^2_b(\R)$, $j\in[m]$, $\varphi\in\Cfo(\R)$, $u\in C_b^2(\R^m)$ and a function $U$ be given by~\eqref{f_U_of_int}. Then $U$ belongs to the domain of the generator $L$ of the Dirichlet form $\e$, that is Friedrich's extension of $(L,\FC)$. Moreover, 
 \begin{equation}\label{f_generator_for_Pp}
 \begin{split}
 LU(g)&=\frac{1}{2}\left[\sum_{i,j=1}^m\partial_i\partial_j u\left(\int_0^1\vec{\alpha}(g(s))ds\right)\cdot\int_0^1\alpha_i'(g(s))\alpha_j'(g(s))ds\right.\\
 &\left.+\sum_{j=1}^m\partial_j u\left(\int_0^1\vec{\alpha}(g(s))ds\right)\cdot\int_0^1\frac{\alpha_j''(g(s))}{m_g(s)}ds\right]\varphi(\|g\|_2^2)\\
 +&\sum_{j=1}^m\partial_ju\left(\int_0^1\vec{\alpha}(g(s))ds\right)\varphi'(\|g\|_2^2)\int_0^1\alpha'(g(s))g(s)ds\\
 +&u\left(\int_0^1\vec{\alpha}(g(s))ds\right)\left[2\varphi''(\|g\|_2^2)\|g\|_2^2+\varphi'(\|g\|_2^2)\cdot\sharp g\right],\quad g\in\St\cap\Li(\xi),
 \end{split}
 \end{equation}
 where $m_g(s)=\leb\{r\in[0,1]:\ g(r)=g(s)\}=\leb g^{-1}(g(s))$, \ $s\in[0,1]$.
\end{proposition}

\begin{proof}
To prove the proposition, it is enough to show that for each $V\in\FC$
 $$
 \e(U,V)=-\langle LU,V\rangle_{\Ll},
 $$
 where $LU$ is defined by~\eqref{f_generator_for_Pp}. The proof of this fact is similar to the proof of Theorem~\ref{theorem_integratin_by_parts}, using the relation $\D U=\pr_{\cdot}\nabla^{L_2}U=\nabla^{L_2}U$. 
\end{proof}

We set
$$
M'_{\alpha}(t):=\langle\alpha,\mu_t\rangle-\langle\alpha,\mu_0\rangle-\int_0^t\Gamma(\mu_s)ds,\quad t\geq 0, 
$$
where $\langle\alpha,\Gamma(\nu)\rangle=\frac{1}{2}\sum_{x\in\supp\nu}\alpha''(x)$,\ \ $\alpha\in\Cfo(\R)$.
Using the martingale problem for $X$ and Proposition~\ref{prop_integr_by_parts_for_mu}, it is easy to obtain the following statement. 

\begin{theorem}\label{theotem_mart_problem_for_M'}
There exists an $\e$-exceptional subset $N$ of $\Li(\xi)$ such that for all $g\in\Li(\xi)\setminus N$, $\alpha\in\Cfo(\R)$ and each $(\F_t)$-stopping time $\tau$ satisfying $\p_{g}\{\tau<\zeta\}=1$ and $\E_gd_{\mathcal W}(\mu_t^{\tau},\leb)^2<\infty$, $t\geq 0$, the process $M_{\alpha}^{\tau}(t)$, $t\geq 0$, is a continuous square integrable $(\F_t)$-martingale under $\p_g$ in $L_2(\xi)$ with the  quadratic variation
$$
\int_0^{t\wedge\tau}\left\langle\left(\alpha'\right)^2,\mu_s\right\rangle ds,
$$ 
where $\mu_t$, $t\geq 0$, is defined by~\eqref{f_mu}, $\mu_t^{\tau}:=\mu_{t\wedge\tau}$ and $M_{\alpha}^{\tau}(t):=M_{\alpha}'(t\wedge\tau)$. 
\end{theorem}

The theorem implies that $\mu_t$, $t\geq 0$, is a martingale solution to equation~\eqref{f_SPDE_for_mu} on $[0,\tau]$.

\begin{corollary}\label{corollary_mart_prop_of_M'}
If $\xi$ is constant on some neighbourhoods of $0$ and $1$, then for $\e$-q.e. $g\in\Li(\xi)$ the process $M_{\alpha}'(t)$, $t\geq 0$, is a continuous square integrable $(\F_t)$-martingale under $\p_g$ in $L_2(\xi)$ with the  quadratic variation
$$
\int_0^t\left\langle\left(\alpha'\right)^2,\mu_s\right\rangle ds.
$$ 
\end{corollary}

\begin{proof}
The corollary follows from Theorem~\ref{theotem_mart_problem_for_M'} and the fact that $\E_g\|X_t\|_2^2<\infty$, $t\geq 0$, for $\e$-q.e. $g\in\Li(\xi)$ (see Corollary~\ref{corollary_M_is_martingale}). 
\end{proof}

Set 
\begin{align}
d_{\mathcal W}(A,B)&=\essinf\{d_{\mathcal W}(\nu_1,\nu_2):\ \nu_1\in A,\ \nu_2\in B\},\\
d_{\mathcal W}(\nu,A)&=\essinf_{\rho\in A}d_{\mathcal W}(\nu,\rho),
\end{align}
for measurable sets $A,B\subset\Pp_2$ and $\nu\in\Pp_2$.

\begin{theorem}
 Let $\xi$ be a strictly increasing function and $\varSigma$ be the push forward of $\Xi$ under the map $\iota$. Then the following statements hold.
 
 $(i)$ For any measurable $A,B\subset\Pp_2$ with $0<\varSigma(A)<\infty$, $0<\varSigma(B)<\infty$ and $A$ or $B$ open we have
 $$
 \lim_{t\to 0}t\ln \int_A\p_{\iota^{-1}\nu}\{\mu_t\in B\}\varSigma(d\nu)=-\frac{d_{\mathcal W}(A,B)^2}{2}.
 $$
  
 \vspace{7pt}
 $(ii)$ Let $A$ be a non empty open  subset of $\Pp_2$ with $\varSigma(A)<\infty$ and let $\Theta$ be any probability measure which is mutually absolutely continuous with respect to $\varSigma$. Then the function $v_t=-t\ln \p_{\iota^{-1}\cdot}\{\mu_t\in A\}$ converges to $\frac{d_{\mathcal W}(\cdot,A)^2}{2}$ in the following senses.
 \begin{enumerate}
  \item[(a)] $v_t\cdot\I_{\{v_t<\infty\}}$ converges to $\frac{d_{\mathcal W}(\cdot,A)^2}{2}\cdot\I_{\{d_{\mathcal W}(\cdot,A)<\infty\}}$ in $\Theta$-probability as $t\to 0$.
  
  \item[(b)] If $F$ is a bounded function on $[0,\infty]$ that is continuous on $[0,\infty)$, then $F(v_t)$ converges to $F\left(\frac{d_{\mathcal W}(\cdot,A)^2}{2}\right)$ in $L_2(\Pp_2,\Theta)$ as $t\to 0$.
 \end{enumerate}
\end{theorem}

\begin{proof}
 The statement follows from Theorems~\ref{theorem_varadhan_formula} and~\ref{theorem_varadhan_formula_spec_case} and the isometry of $\Li(\xi)=\Li$ and $\Pp_2$. 
\end{proof}

\appendix

\section{Appendix}

\subsection{\texorpdfstring{$\Li(\xi)$}{Li(xi)}-functions}

 Let $\xi$ be a bounded 
 function from $D^{\uparrow}$. Recall that $\Li(\xi)$ denote the set of functions from $\Li$ that are $\sigma^{\star}(\xi)$-measurable.
 
 \begin{remark}
  \begin{enumerate}
   \item[(i)] The space $\Li(\xi)$ is closed in $\Li$.
   
   \item[(ii)] Let $f\in\Li(\xi)$ and $g$ be its modification from $D^{\uparrow}$, then $g$ is $\sigma^{\star}(\xi)$-measurable.
  \end{enumerate}

 \end{remark}
 
 In this section, we will give a useful description of each function $g\in\Li(\xi)$ using its version from $D^{\uparrow}$ denoted also by $g$.
 
 \begin{proposition}\label{prop_xi_measurability}
  A function $g\in\Li$ belongs to $\Li(\xi)$ if and only if for all $a<b$ from $[0,1]$ the equality $\xi(a)=\xi(b)$ implies $g(a)=g(b-)$. 
 \end{proposition}

 \begin{proof}
 Let $g\in\Li(\xi)$ and $\xi(a)=\xi(b)$ for some $a<b$ and $f$ is $\sigma(\xi)$ measurable with $g=f$ a.e. We note that the sets
 $$
 \pi_r=\xi^{-1}(\{r\})=\{s\in[0,1]:\xi(s)=r\},
 $$
 are the smallest in $\sigma(\xi)$, i.e. for any non empty set $A\in\sigma(\xi)$ satisfying $A\subseteq\pi_r$ we have $A=\pi_r$. Consequently, the set 
 $$
 B=\{s\in[0,1]:\ f(a)=f(s)\}\cap\pi_{\xi(a)}
 $$ 
 coincides with $\pi_{\xi(a)}$. We next remark that $[a,b]\subseteq \pi_{\xi(a)}=B$, since $\xi$ is non decreasing and $\xi(a)=\xi(b)$. Therefore, $f(a)=f(s)$ for all $s\in [a,b]$. Thus, the equality $f=g$ a.e. yields $g(a)=g(a+)=g(b-)$.
 
 To prove the sufficiency, we first show that a function $f$ is $\sigma(\xi)$ measurable, if $f$ is Borel measurable and 
 \begin{equation}\label{f_prop_of_measurability}
 \xi(a)=\xi(b)\ \ \mbox{implies}\ \  f(a)=f(b)\ \ \mbox{for all}\ \  a,b\in[0,1].
 \end{equation}
 Let us define the function $\eta[\xi(0),\xi(1)]\to[0,1]$, that will play a role of the inverse function for $\xi$, as follows
 $$
 \eta(r)=\min\{s\in[0,1]:\ \xi(s)\geq r\},\quad r\in[\xi(0),\xi(1)].
 $$
 Then it is easy to see that $\eta$ satisfies the following properties
 \begin{enumerate}
  \item[a)] $\eta$ is a non decreasing left-continuous function;
  
  \item[b)] $\eta(\xi(s))=\widetilde{s}$, where $\widetilde{s}=\min\{\pi_{\xi(s)}\}$.
 \end{enumerate}
 
 Using these properties and setting $\phi(r)=f(\eta(r))$, $r\in[\xi(0),\xi(1)]$, we can easily see that $\phi$ is a Borel function and
 $$
 \phi(\xi(s))=f(\eta(\xi(s)))=f(\widetilde{s})=f(s),\quad s\in[0,1].
 $$
 Thus, $f$ is $\sigma(\xi)$-measurable, as a compositions of Borel function with $\xi$.
 
 Let for all $a<b$ the equality $\xi(a)=\xi(b)$ implies $g(a)=g(b-)$. We are going to find a function $f$ that satisfies~\eqref{f_prop_of_measurability} and coincides with $g$ a.e. Denote the set of all discontinuous points of $g$ by $D_g$ that is at most countable, since $g$ is non decreasing. Next, for all $b\in D_g$ we note that $b$ satisfies only one of the following properties
 \begin{itemize}
  \item $\xi(a)\neq \xi(b)$ for all $a\neq b$;
  
  \item there exists $a<b$ such that $\xi(a)=\xi(b)$ and, consequently, $g(a)=g(b-)$;
  
  \item there exists $c>b$ such that $\xi(b)=\xi(c)$ and, consequently, $g(b)=g(c-)$.
 \end{itemize}
 Indeed, if there exist both $a$ and $c$ such that $a<b<c$ and $\xi(a)=\xi(b)=\xi(c)$ then $g(a)=g(c-)$. But it contradicts the assumption that $b$ is a discontinuous point of $g$. 

 We define 
 \begin{align}
  f(s)=\begin{cases}
        g(s),&\mbox{if}\ \ s\in[0,1]\setminus D_g,\\
        g(s),&\mbox{if}\ \ s\in D_g\ \ \mbox{and}\ \ \xi(a)=\xi(s)\ \ \mbox{for some}\ \ a<s,\\
        g(s-),&\mbox{if}\ \ s\in D_g\ \ \mbox{and}\ \ \xi(s)=\xi(c)\ \ \mbox{for some}\ \ c>s.
       \end{cases}
 \end{align}
 Then $f$ is a well-defined non decreasing function and, consequently, Borel measurable. Moreover, it is easily seen that $f$ satisfies~\eqref{f_prop_of_measurability}. So, $f$ is $\sigma(\xi)$-measurable. Since $D_g$ is at most countable and $\{s:\ g(s)\neq f(s)\}\subseteq D_g$, we have that $f=g$ a.e. Thus, $g$ is $\sigma^{\star}(\xi)$-measurable, that completes the proof of the proposition. 
 \end{proof}
 
 \subsection{Multivariate Bernstein polynomials}
 
 In this section, we give a slight modification of the result obtained in~\cite{Veretennikov:2016} about uniform approximation of a function and its partial derivatives by Bernstein polynomials.
 
 For a function $f:[0,1]^k\to\R$ we define the Bernstein polynomials on $[0,1]^k$ as follows 
 \begin{align}
  B_n(f;x)&=\sum_{j_1,\ldots,j_k=0}^nf\left(\frac{j_1}{n},\ldots,\frac{j_k}{n}\right)C_n^{j_1}\ldots C_n^{j_k}\\
  &\cdot x_1^{j_1}(1-x_1)^{n-j_1}\ldots x_k^{j_k}(1-x_k)^{n-j_k},
 \end{align}
 where $C_n^j=\frac{n!}{j!(n-j)!}$,\ \  $j\in[n]\cup\{0\}$.

  \begin{proposition}
   If $f\in C^1(\R^k)$, then 
   \begin{enumerate}
    \item[(i)] $\{B_n(f;\cdot)\}_{n\geq 1}$ uniformly converges to $f$ on $[0,1]^k$;
    
    \item[(ii)] $\{\partial_iB_n(f;\cdot)\}_{n\geq 1}$ uniformly converges to $\partial_if$ on $[0,1]^k$ for all $i\in[k]$.
   \end{enumerate}
  \end{proposition}
  \begin{proof}
   The statement is a partial case 
   of Theorem~4~\cite{Veretennikov:2016}. 
  \end{proof}

  Next we would like to have a sequence of polynomials that approximate a function $f$ on $[-M,M]^k$. We set for a fixed $M>0$
  \begin{align}\label{f_polynomial_P}
  \begin{split}
   f_M(x)&=f(2Mx-M),\\
   P_n^M(f;x)&=B_n\left(f_M;\frac{x}{2M}+\frac{1}{2}\right)-B_n\left(f_M;\frac{1}{2}\right).
   \end{split}
  \end{align}
  
  We note that $P_n^M(f;0)=0$. This property is important for us, since in this case the composition $P_n^M(f;U)$ belongs to $\FC$ for $U_i\in\FC$, $i\in[k]$.
  
  The following proposition is a trivial consequence of the previous proposition.
  
  \begin{lemma}\label{lemma_converg_of_polynomials}
   Let $f\in C^1(\R^k)$ and $f(0)=0$. Then
   \begin{enumerate}
    \item[(i)] $\{P_n^M(f;\cdot)\}_{n\geq 1}$ uniformly converges to $f$ on $[-M,M]^k$;
    
    \item[(ii)] $\{\partial_iP_n^M(f;\cdot)\}_{n\geq 1}$ uniformly converges to $\partial_if$ on $[-M,M]^k$ for all $i\in[k]$.
   \end{enumerate}

  \end{lemma}

 \subsection{Proof of auxiliary statements}
 \subsubsection{Proof of Lemma~\ref{lem_proj_of_step_function}}
 By Remark~\ref{rem_proj_and_ecpect}~(iii), $\pr_gh$ belongs to $\Li$. Thus, we need only to show that it has a modification that takes a finite number of values. Consequently, using the linearity of $\pr_g$ and Remark~\ref{rem_repres_of_step_functions}, it is enough to prove that for any $H:=[a,b)\subset[0,1]$, $\pr_g\I_H$ has a modification that takes at most three values. 

 We set $D_n=\left\{\frac{k}{2^n},\ k\in\Z\right\}$, $\s_n=\sigma\{[a,b):\ a<b,\ a,b\in D_n\}$ and $\F_n=g^{-1}(\s_n)$. Let us note that $\{\F_n,\ n\in\N\}$, is increasing, since $\{\s_n,\ n\in\N\}$ increases. Moreover, it is clear that 
 $$
 \sigma(g)=\bigvee_{n=1}^{\infty}\F_n=\sigma\left(\bigcup_{n=1}^{\infty}\F_n\right).
 $$
 By Levi's theorem (see, e.g.,~\cite[Theorem~1.5]{Liptser:2001}),
 \begin{equation}\label{f_conv_of_cond_expect}
 \E(\I_H|\F_n)\to\E\left(\I_H\left|\bigvee_{n=1}^{\infty}\F_n\right.\right)\quad \mbox{a.e., \ as}\ n\to\infty,
 \end{equation}
 where $\E$ denotes the expectation on the probability space $([0,1],\B([0,1]),\leb)$. 
 Since each element of $\F_n$ can be written as a finite or a countable union of disjoint sets $G_{k,n}=g^{-1}\left(\left[\frac{k}{2^n},\frac{k+1}{2^n}\right)\right)$, $k\in\Z$, we obtain
 $$
 \E(\I_H|\F_n)=\sum_{k\in\Z}\frac{\I_{G_{k,n}}}{\leb(G_{k,n})}\E\I_{H\cap G_{k,n}}.
 $$
 Next, by monotonicity of $g$, the set $H$ can be covered by a finite number of $G_{k,n}$, i.e there exist integer numbers $p_1<p_2$ such that
 \begin{itemize}
  \item $\widetilde{H}:=\bigcup_{k=p_1+1}^{p_2-1}G_{k,n}\subseteq H=[a,b)$;
  
  \item $a\in G_{p_1,n}$, $b\in G_{p_2,n}$;
  
  \item for each $k<p_1$ or $k>p_2$, $G_{k,n}\cap H=\emptyset$.
 \end{itemize}
 Thus,
 $$
 \E(\I_H|\F_n)=\frac{\I_{G_{p_1,n}}}{\leb(G_{p_1,n})}\E\I_{H\cap G_{p_1,n}}+\frac{\I_{G_{p_2,n}}}{\leb(G_{p_2,n})}\E\I_{H\cap G_{p_2,n}}+\frac{\I_{\widetilde{H}}}{\leb(\widetilde{H})}\E\I_{\widetilde{H}}.
 $$
 Hence $\E(\I_H|\F_n)$ takes at most three values. By~\eqref{f_conv_of_cond_expect} and Remark~\ref{rem_proj_and_ecpect}, $\pr_g\I_H$ also takes at most three values. This completes the proof of the lemma.

 \subsubsection{Proof of Proposition~\ref{prop_int_alpha_to_dom}}
 Note that the sequence of $\sigma$-algebras 
 $$
 \s_n=\sigma\left(\pi_i^n:=\left[\frac{i-1}{2^n},\frac{i}{2^n}\right),\ i\in\left[2^n\right]\right),\quad n\in\N,
 $$
 increases to $\B([0,1])$, that is, $\sigma\left(\bigcup_{n=1}^\infty\s_n\right)=\B([0,1])$. Considering functions from $\Li(\xi)$ as random elements on the probability space $([0,1],\B([0,1]),\leb)$ and using the Levy theorem (see, e.g.,~\cite[Theorem~1.5]{Liptser:2001}), for each $g\in\Li(\xi)$
 $$
 g_n:=\E(g|\s_n)=\sum_{i=1}^{2^n}\langle g,h_i^n\rangle\I_{\pi_i^n}\to g\quad \mbox{a.s. \ as}\ \ n\to\infty,
 $$
 where $h_i^n=2^n\I_{\pi_i^n}$. Therefore, by the dominated convergence theorem,
 $$
 \int_0^1\alpha(g_n(s))ds=\sum_{i=1}^{2^n}\alpha(\langle g,h_i^n\rangle)\frac{1}{2^n}\to \int_0^1\alpha(g(s))ds\quad \mbox{as}\ \ n\to\infty.
 $$
 
 We next define 
 $$
 U_n(g)=\int_0^1\alpha(g_n(s))ds\cdot\varphi(\|g\|_2^2),\quad g\in\Li(\xi),
 $$
 and note that $U_n\in\FC$. Moreover, for all $g\in\Li(\xi)$
 \begin{align}
 \D U_n(g)&=\frac{1}{2^n}\sum_{i=1}^{2^n}\alpha'(\langle g,h_i^n\rangle)\pr_gh_i^n\varphi(\|g\|_2^2)+2\int_0^1\alpha(g_n(s))ds\cdot\varphi'(\|g\|_2^2)g\\
 &=\pr_g\alpha'(g_n)\varphi(\|g\|_2^2)+2\int_0^1\alpha(g_n(s))ds\cdot\varphi'(\|g\|_2^2)g.
 \end{align}
 By the dominated convergence theorem and Remark~\ref{rem_proj_and_ecpect}~(ii), 
 $$
 \pr_g\alpha'(g_n)=\E(\alpha'(g_n)|\sigma^{\star}(g))\to\E(\alpha'(g)|\sigma^{\star}(g))=\alpha'(g)\quad\mbox{a.s. \ as}\ \ n\to\infty.
 $$
 Thus, using the dominated convergence theorem again, we have
 $$
 U_n\to U\quad\mbox{and}\quad \|\D U_n-\D U\|_2\to 0\quad\mbox{in}\ \ \Ll\ \ \mbox{as}\ \ n\to\infty,
 $$
 where $U(g)=\int_0^1\alpha(g(s))ds\cdot\varphi(\|g\|_2^2)$ and $\D U(g)=\alpha'(g)\varphi(\|g\|_2^2)+2\int_0^1\alpha(g(s))ds\cdot\varphi'(\|g\|_2^2)g$,\ \  $g\in\Li(\xi)$. The proposition is proved.
 
\subsection*{Acknowledgements}
The research of the first author was partly supported by Alexander von Humboldt Foundation and partly supported by the Deutsche
Forschungsgemeinschaft (DFG, German Research Foundation) – SFB 1283/2 2021 – 317210226.

\nocite{MR2842966,Arratia:1979,MR2520126,MR2094432,Le_Jan:2004,Dawson:1993,MR1612725,MR2994690,MR1915445}


\providecommand{\bysame}{\leavevmode\hbox to3em{\hrulefill}\thinspace}
\providecommand{\MR}{\relax\ifhmode\unskip\space\fi MR }
\providecommand{\MRhref}[2]{%
  \href{http://www.ams.org/mathscinet-getitem?mr=#1}{#2}
}
\providecommand{\href}[2]{#2}

\end{document}